\pdfoutput=1
\documentclass[pdflatex,sn-mathphys-num]{sn-jnl}
\usepackage{graphicx} %

\usepackage{amsthm}
\usepackage{amsfonts}
\usepackage{amssymb}
\usepackage{amsmath}
\usepackage{xspace}
\usepackage{mathrsfs}
\usepackage{tikz}
\usepackage{tikz-cd}
\hypersetup{hidelinks}

\usepackage{booktabs}
\usepackage{siunitx}
\newcommand{\?}{\mskip1.5mu}
\newcommand{\sing}[1]{\left\{\?#1\?\right\}}

\newcommand{\N}{\mathbf N}
\newcommand{\R}{\mathbf R}

\newcommand{\E}{\mathscr E}
\newcommand{\I}{\mathscr I}

\newcommand{\Dy}{\mathscr D}
\newcommand{\Stair}{\mathscr S}
\newcommand{\downset}{\mathord\downarrow}
\newcommand{\upset}{\mathord\uparrow}

\renewcommand{\emptyset}{\varnothing}

\newcommand{\sqleq}{\sqsubseteq}
\newcommand{\sqless}{\sqsubset}
\newcommand{\sqgeq}{\sqsupseteq}

\newcommand{\OR}{\vee}
\newcommand{\AND}{\wedge}

\def\..{\,\mathpunct{\ldotp\ldotp}} %

\newcommand{\dw}[1]{\mathsf{#1}}
\newcommand{\CCB}{Clarke--Cormack--Burkowski\xspace}

\newcommand{\area}{\mathsf{area}}
\newcommand{\rank}{\mathsf{rank}}
\newcommand{\sweepdim}{\mathsf{sweepdim}}
\renewcommand{\dim}{\mathsf{dim}}
\newcommand{\width}{\mathsf{width}}

\theoremstyle{thmstyleone}
\newtheorem{theorem}{Theorem}
\newtheorem{proposition}{Proposition}
\newtheorem{lemma}{Lemma}
\newtheorem{corollary}{Corollary}
\theoremstyle{thmstylethree}
\newtheorem{definition}{Definition}
\newtheorem{remark}{Remark}
\renewcommand\epsilon\varepsilon

\title[Antichains of intervals, Dyck paths, and sweep dimension]{The lattice of antichains of intervals: Dyck paths, level sets, and sweep dimension}

\author*[1]{\fnm{Paolo} \sur{Boldi}}\email{paolo.boldi@unimi.it}
\author[1]{\fnm{Sebastiano} \sur{Vigna}}\email{sebastiano.vigna@unimi.it}
\affil*[1]{\orgdiv{Dipartimento di Informatica}, \orgname{Università degli Studi di Milano}, \orgaddress{\city{Milan}, \country{Italy}}}

\abstract{The \CCB lattice $\E_n$ is the lattice of antichains of intervals of the chain $\{\?0,1,\dots,n-1\?\}$, the algebraic structure underlying minimal-interval semantics, a query semantics used in information retrieval. We give an explicit, elementary anti-isomorphism between $\E_n$ minus its top element and the Stanley lattice of Dyck words of semilength $n+1$ (an anti-isomorphism whose existence is implicit in work of Ferrari and Munarini), under which the rank of an antichain is $\binom{n+1}{2}$ minus the area of the associated Dyck word. Hence the rank generating function of $\E_n\setminus\{\top\}$ is the reverse of the Carlitz--Riordan $q$-Catalan polynomial. Using a local limit theorem for the area of random-walk excursions, we confirm unconditionally an asymptotic formula of Tak\'acs: the size $\lambda_n$ of the largest level set of $\E_n$ satisfies $\lambda_n\sim\kappa\,4^{n+1}/(n+1)^3$, where $\kappa\approx1.0807$. Finally, for the antichains of intervals of an arbitrary finite poset $P$ we study the \emph{sweep dimension}, the minimum number of linear orders needed to separate every interval from every element outside it. It equals the Dushnik--Miller dimension of the lattice of convex subsets of $P$, it is the least $k$ for which $P$ maps into $\R^k$ so that every interval is the preimage of an upper orthant, and it lies between $\max(2,w)$ and $2w$, where $w$ is the width of $P$, both bounds being tight. Antichains of intervals embed, preserving joins, into $k$-dimensional staircase regions.}

\keywords{antichains of intervals, minimal-interval semantics, Dyck paths, Stanley lattice, Carlitz--Riordan $q$-Catalan numbers, order dimension}

\pacs[MSC Classification]{06A07, 06D05, 05A15, 05A16, 68P20}

\begin{document}

\maketitle

\section{Introduction}

In modern information-retrieval systems, such as web search engines, the fact that a document (seen as a sequence of words) satisfies a query is not enough to determine its relevance. In fact, the user is often interested in the \emph{regions} of the document that satisfy the query, and in particular in the \emph{minimal} regions that satisfy it. For instance, if a user searches for ``\texttt{(A AND B) OR (C AND D)}'', it is not enough to know that a document contains both $A$ and $B$ (or both $C$ and $D$): the user wants to know where these terms occur, and in particular whether they occur close to each other.

This idea is formally captured by \emph{minimal-interval
semantics}~\cite{CCBASTSFI}, which uses \emph{antichains} (w.r.t.~inclusion) of
\emph{intervals} (i.e., textual passages in the document) to represent the
answer to a query. Each interval in the text of a document is a \emph{witness}
of the satisfiability of the query, that is, it represents a region of the
document that satisfies the query; only minimal witnesses are kept, since they
are the most informative ones. 

This semantics belongs to a larger family of set-based algebraic query languages for text
indexing systems that record, besides the positions of words, the positions of
\emph{text regions}; Consens and Milo~\cite{CoMAQTR} list several of them,
including the algebra of the Pat system and the one of Clarke, Cormack and
Burkowski~\cite{CCBASTSFI}, and study the expressive power and the
optimization of the Pat algebra, which, in their words, includes most of the
operators considered in the other proposals. In the approach of Clarke, Cormack
and Burkowski~\cite{CCBASTSFI}, the algebraic laws needed by the query
operators (e.g., distributivity) are proved from scratch; minimal intervals can
also be used to estimate the relevance of a document~\cite{ClCSSRR}.
In~\cite{BoVLAFI} we showed that these laws are instances of general
properties of \emph{antichain completions}: the antichains of intervals of a
document form a lattice, which we call the \CCB lattice, in which the Boolean
operators of queries become lattice operations; there we also
characterized in elementary terms the relative pseudo-complement and other
operators of interest in information retrieval, and gave a closed formula for
the rank of an element; efficient, optimally lazy algorithms for the operators
were described in~\cite{BoVEOLAMIS}. More recently, Clarke's \emph{annotative
indexing}~\cite{ClaAI} has extended minimal-interval semantics in two
directions: arbitrary intervals, and not only single positions, can be
indexed, and every interval carries a value; for each feature, the set of
annotations is still required to follow minimal-interval semantics, that is,
to be an element of a \CCB lattice.

\smallskip
In this paper, we give an explicit anti-isomorphism between the lattice of antichains of intervals underlying minimal-interval semantics, for a total order, and the Stanley lattice of Dyck paths, which yields the exact asymptotics of its largest level set. For an arbitrary finite poset, instead, we show that it can be embedded into a lattice of higher-dimensional staircase regions whose dimension is a parameter of the poset, the sweep dimension, which turns out to be the Dushnik--Miller dimension of the lattice of convex subsets of the poset; we characterize it and bound it in terms of the width.

\paragraph{What is, and what is not, new.} Section~\ref{sec:related} discusses the literature in detail; in short: the existence of an anti-isomorphism between the lattice (deprived of its top element) and the Stanley lattice is implicit in the work of Ferrari and Munarini~\cite{FeMLPRTV} through Birkhoff's representation theorem, and what we add is an explicit, elementary anti-isomorphism with complete proofs and the dictionary between rank and area; the asymptotic formula for the largest level set, with its constant, is due to Tak\'acs~\cite{TakBEVA}, who proved it under an assumption of unimodality and of a local limit theorem, and what we add is to make it unconditional through the local limit theorem of Denisov, Kolb and Wachtel~\cite{DKWLAARWE} and to transfer it to $\E_n$; the sweep dimension of a poset turns out to be the Dushnik--Miller dimension of its lattice of convex subsets (Proposition~\ref{prop:dimsweep}), and its bounds in terms of the width can also be recovered from the work of Edelman and Saks~\cite{EdSCRCDCG} on convex geometries, whereas the orthant and staircase-region descriptions, and the direct proofs yielding explicit sweep systems, are new.

\paragraph{Organization of the paper.} After setting up the notation and conventions (Section~\ref{sec:notation}), in Section~\ref{sec:ccb} we focus on the \CCB lattice $\E_n$; more precisely, we prove that $\E_n$ minus its top element
is anti-isomorphic to the Stanley lattice of Dyck words of semilength $n+1$,
through the intermediate representation by staircase functions. We also prove that the
rank of an antichain is $\binom{n+1}2$ minus the area of the corresponding Dyck
word; as a consequence, the level sets of $\E_n$ are counted by the
coefficients of the Carlitz--Riordan $q$-Catalan polynomial and by partitions
contained in a staircase; this observation is used in Section~\ref{sec:asymptotics}, where we derive the
asymptotics of the largest level set from a local limit theorem for the area of
random-walk excursions. Section~\ref{sec:ext} extends the picture to the
antichain lattice of intervals of an arbitrary finite poset: we introduce the
sweep dimension, characterize it through representations of the poset in
$\R^k$ in which intervals are preimages of orthants, obtain the join-preserving
embedding of antichains of intervals into $k$-dimensional staircase regions,
prove that the sweep dimension is the order dimension of the poset of intervals,
and bound it in terms of the width, showing that the bounds are tight. Section~\ref{sec:related} compares our results with the literature on lattices of Dyck paths, on the unimodality of the Carlitz--Riordan polynomials and on the dimension of posets of convex sets, and Section~\ref{sec:conclusions} wraps up the paper with some open problems.

\section{Definitions and notations}
\label{sec:notation}

We briefly recall the setting and notation of~\cite{BoVLAFI}, to which we refer for proofs and for the general theory. 

A \emph{partial order} $\sqleq$ on a set $P$ is a binary relation that is reflexive, antisymmetric, and transitive; the pair $(P,\sqleq)$ is called a \emph{poset} (partially ordered set). We say that $\sqleq$ is \emph{linear} (or total) if for any $x,y \in P$,  $x \sqleq y$ or $y \sqleq x$.

\smallskip
A \emph{chain} of $(P,\sqleq)$ is a set of pairwise comparable elements of $P$. An \emph{antichain} of $(P,\sqleq)$ is a set of pairwise incomparable elements of $P$.  The cardinality of the largest antichain of a finite poset $(P,\sqleq)$ is called its \emph{width}, and it is denoted by $\width(P,\sqleq)$: by Dilworth's theorem~\cite{DilDTPOS}, the minimum number of chains in a partition of a finite poset $P$ (a \emph{chain decomposition}) equals $\width(P,\sqleq)$.

An \emph{interval} of a poset $(P,\sqleq)$ is a set $I \subseteq P$ such that for any $x,y,z \in P$, if $x \sqleq y \sqleq z$ and $x,z \in I$ then $y \in I$. The set of intervals of $(P,\sqleq)$ is denoted by $\I(P,\sqleq)$. 

Intervals are ordered by reverse inclusion ($\supseteq$), and the set of antichains of intervals (i.e., sets of pairwise inclusion-incomparable intervals) is denoted by $\E(P,\sqleq)$. The order on $\E(P,\sqleq)$ is defined as follows:
\[
  A \leq B \iff \forall I \in A, \exists J \in B \text{ such that } J \subseteq I.
\]
The intended meaning of this order is that $A \leq B$ if every interval in $A$ contains an interval in $B$, that is, if $B$ is ``more informative'' than $A$ (for every witness of $A$, there is an even more informative witness in $B$, and $B$ can also contain additional intervals).

In this paper we shall \emph{only consider finite posets} $(P,\sqleq)$: under this condition\footnote{In fact, under much more general conditions, but in this paper we will restrict ourselves to the finite case.} we have~\cite{BoVLAFI}:
\begin{theorem}
  If $(P,\sqleq)$ is finite, then $\E(P,\sqleq)$ is a strongly coatomic, superalgebraic, completely distributive lattice.
\end{theorem}

\section{The \CCB lattice}\label{sec:ccb}

If $n$ is a natural number, we let $[n]=\{\?0,1,\dots,n-1\?\}$ with the usual order. This is a total order, and will play the role of $(P,\sqleq)$ in our previous definitions. We let $\I_n$ be a shortcut for $\I([n],\leq)$, and $\E_n=\E([n],\leq)$ be the lattice of interval antichains, called the \emph{\CCB lattice}.

Throughout the paper, a chain with $n$ elements is related to Dyck words of semilength $n+1$ (see Theorem~\ref{thm:delta}), and quantities indexed by the semilength $m=n+1$ appear shifted by one with respect to those indexed by $n$.

We use the Hoare--Ramshaw notation $[\ell\..r]=\{\?x\mid \ell\leq
x\leq r\?\}$ for nonempty intervals of $[n]$. 

We recall some fundamental properties of $\E_n$ (see~\cite{BoVLAFI}); in particular, being an instance of the previous theorem, $\E_n$ is a complete lattice, and we use $\OR$ and $\AND$ to denote its join and meet operations, respectively.

\begin{theorem}
  For every $n \geq 1$, the following properties of the \CCB lattice $\E_n$ hold:
  \begin{enumerate} 
    \item The minimum $\bot$ of $\E_n$ is the empty antichain $\emptyset$.
    \item The maximum $\top$ of $\E_n$ is the antichain $\sing{\emptyset}$ containing the empty interval.
    \item Every element of $\E_n$, except $\top$, contains only nonempty intervals.
    \item The only atom of $\E_n$ is the antichain $\sing{[0\..n-1]}$, whose only interval is the whole $[n]$.
    \item The only coatom of $\E_n$ is the antichain $\sing{[x\..x] \mid x \in [n]}$, containing all the size-1 intervals.
  \end{enumerate}
\end{theorem}

\begin{figure}
	\centering
	\includegraphics[scale=.4]{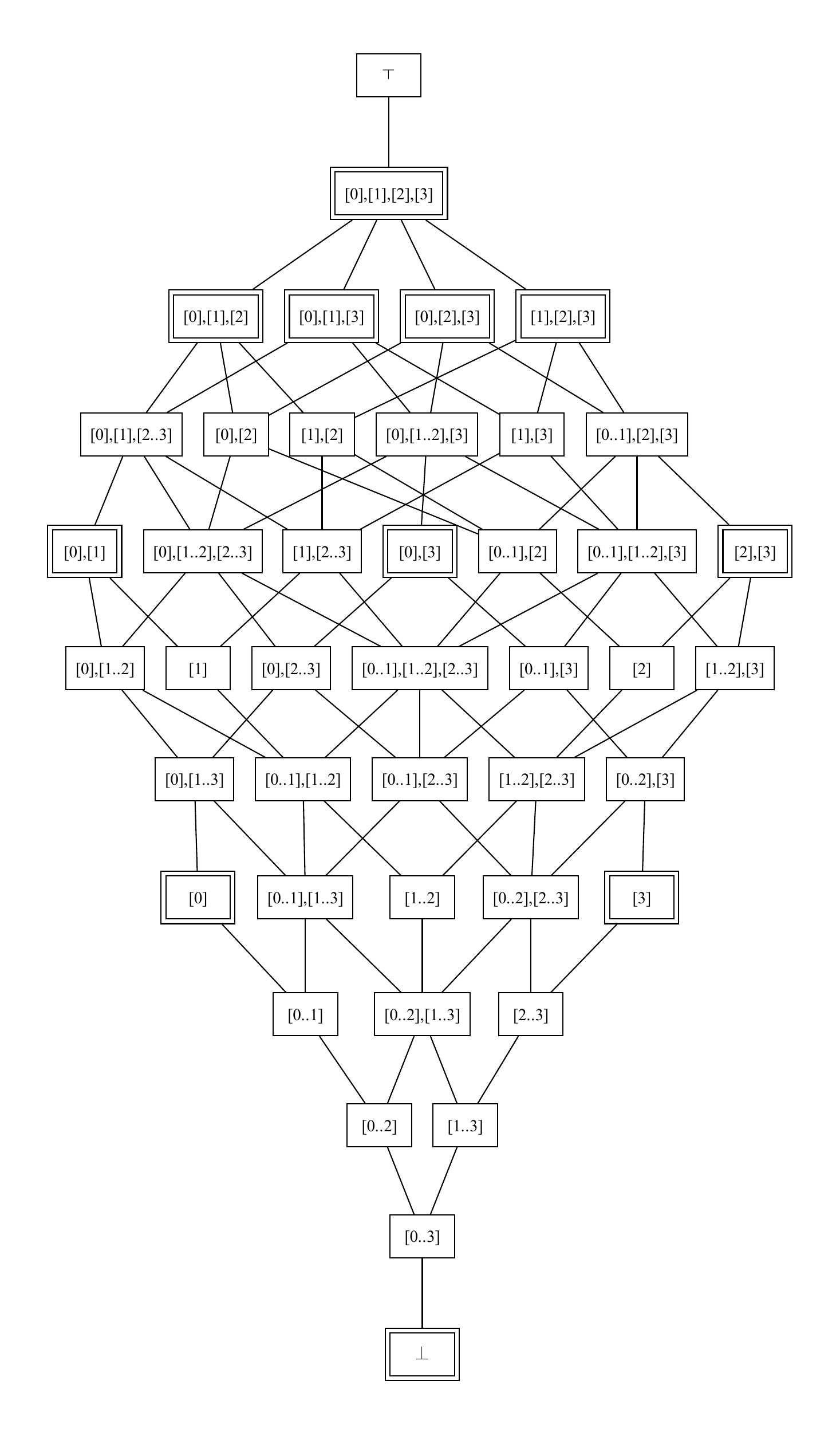}
	\caption{\label{fig:e4}The lattice $\E_4$ (represented by its Hasse diagram). To improve readability, we omitted the outer brackets of the antichains. For instance, the element $\sing{[0\..1],[2\..3]}$ is represented as $[0\..1],[2\..3]$. Elements with double frame are $\AND$-irreducible.}
\end{figure}

\subsection{Dyck words, dominance order and staircase functions}

A \emph{Dyck word} of semilength $m$ is a word $w$ of length $2m$ over the alphabet $\{\?\dw U,\dw D\?\}$ containing $m$ occurrences of each letter and such that no prefix contains more $\dw D$'s than $\dw U$'s. We denote by $\Dy_m$ the set of Dyck words of semilength $m$. The \emph{height function}
of $w\in\Dy_m$ is the function $h_w:[2m+1]\to\N$ where $h_w(t)$
is the number of $\dw U$'s minus the number of $\dw D$'s in the prefix of $w$ of length $t$; drawing $\dw U$ as a rise and $\dw D$ as a fall, $h_w$ is just
the profile of the path associated with $w$. Clearly $w\mapsto h_w$ is injective, $h_w(0)=h_w(2m)=0$, $h_w(t)\geq0$, and $|h_w(t+1)-h_w(t)|=1$ for all $t$; conversely, every function with these properties is the height function of a (unique) Dyck word.

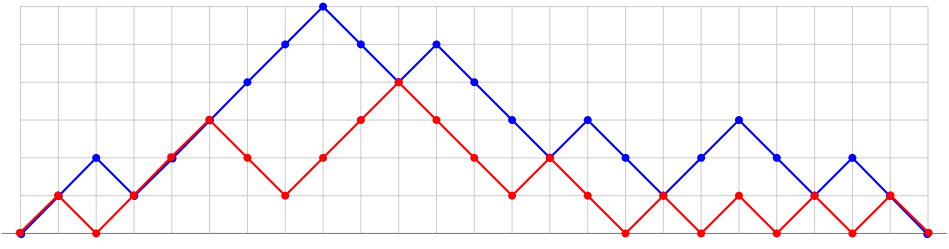
\begin{figure}
  \centering
  \begin{tikzpicture}[scale=.5]
  \draw[very thin, gray!40] (0,0) grid (24,6);
  \draw[gray] (-.5,0) -- (24.5,0);
  \draw[thick, blue] (0.018,-0.018) -- (1.009,0.991) -- (2.000,2.000) -- (3.009,0.991) -- (4.018,1.982) -- (5.009,2.991) -- (6.000,4.000) -- (7.000,5.000) -- (8.000,6.000) -- (9.000,5.000) -- (10.000,4.000) -- (11.000,5.000) -- (12.000,4.000) -- (13.000,3.000) -- (14.000,2.000) -- (15.000,3.000) -- (16.000,2.000) -- (17.000,1.000) -- (18.000,2.000) -- (19.000,3.000) -- (20.000,2.000) -- (21.000,1.000) -- (22.000,2.000) -- (22.991,0.991) -- (23.982,-0.018);
  \fill[blue] (0.018,-0.018) circle (3pt);
  \fill[blue] (1.009,0.991) circle (3pt);
  \fill[blue] (2.000,2.000) circle (3pt);
  \fill[blue] (3.009,0.991) circle (3pt);
  \fill[blue] (4.018,1.982) circle (3pt);
  \fill[blue] (5.009,2.991) circle (3pt);
  \fill[blue] (6.000,4.000) circle (3pt);
  \fill[blue] (7.000,5.000) circle (3pt);
  \fill[blue] (8.000,6.000) circle (3pt);
  \fill[blue] (9.000,5.000) circle (3pt);
  \fill[blue] (10.000,4.000) circle (3pt);
  \fill[blue] (11.000,5.000) circle (3pt);
  \fill[blue] (12.000,4.000) circle (3pt);
  \fill[blue] (13.000,3.000) circle (3pt);
  \fill[blue] (14.000,2.000) circle (3pt);
  \fill[blue] (15.000,3.000) circle (3pt);
  \fill[blue] (16.000,2.000) circle (3pt);
  \fill[blue] (17.000,1.000) circle (3pt);
  \fill[blue] (18.000,2.000) circle (3pt);
  \fill[blue] (19.000,3.000) circle (3pt);
  \fill[blue] (20.000,2.000) circle (3pt);
  \fill[blue] (21.000,1.000) circle (3pt);
  \fill[blue] (22.000,2.000) circle (3pt);
  \fill[blue] (22.991,0.991) circle (3pt);
  \fill[blue] (23.982,-0.018) circle (3pt);
  \draw[thick, red] (-0.018,0.018) -- (0.991,1.009) -- (2.000,0.000) -- (2.991,1.009) -- (3.982,2.018) -- (4.991,3.009) -- (6.000,2.000) -- (7.000,1.000) -- (8.000,2.000) -- (9.000,3.000) -- (10.000,4.000) -- (11.000,3.000) -- (12.000,2.000) -- (13.000,1.000) -- (14.000,2.000) -- (15.000,1.000) -- (16.000,0.000) -- (17.000,1.000) -- (18.000,0.000) -- (19.000,1.000) -- (20.000,0.000) -- (21.000,1.000) -- (22.000,0.000) -- (23.009,1.009) -- (24.018,0.018);
  \fill[red] (-0.018,0.018) circle (3pt);
  \fill[red] (0.991,1.009) circle (3pt);
  \fill[red] (2.000,0.000) circle (3pt);
  \fill[red] (2.991,1.009) circle (3pt);
  \fill[red] (3.982,2.018) circle (3pt);
  \fill[red] (4.991,3.009) circle (3pt);
  \fill[red] (6.000,2.000) circle (3pt);
  \fill[red] (7.000,1.000) circle (3pt);
  \fill[red] (8.000,2.000) circle (3pt);
  \fill[red] (9.000,3.000) circle (3pt);
  \fill[red] (10.000,4.000) circle (3pt);
  \fill[red] (11.000,3.000) circle (3pt);
  \fill[red] (12.000,2.000) circle (3pt);
  \fill[red] (13.000,1.000) circle (3pt);
  \fill[red] (14.000,2.000) circle (3pt);
  \fill[red] (15.000,1.000) circle (3pt);
  \fill[red] (16.000,0.000) circle (3pt);
  \fill[red] (17.000,1.000) circle (3pt);
  \fill[red] (18.000,0.000) circle (3pt);
  \fill[red] (19.000,1.000) circle (3pt);
  \fill[red] (20.000,0.000) circle (3pt);
  \fill[red] (21.000,1.000) circle (3pt);
  \fill[red] (22.000,0.000) circle (3pt);
  \fill[red] (23.009,1.009) circle (3pt);
  \fill[red] (24.018,0.018) circle (3pt);
\end{tikzpicture}
  \caption{\label{fig:dyckword}The height function of the Dyck words $w$ ($\dw{UUDUUUUUDDUDDDUDDUUDDUDD}$, blue) and $v$ ($\dw{UDUUUDDUUUDDDUDDUDUDUDUD}$,  red). Here $v \preceq w$ according to the dominance order: the blue path lies weakly above the red one.}
\end{figure}

The \emph{dominance order} on $\Dy_m$ is defined by
\[
  v\preceq w\qquad\text{iff}\qquad h_v(t)\leq h_w(t)\text{ for all $t$},
\]
that is, $v\preceq w$ iff the path of $v$ lies weakly below the path of $w$. The poset $(\Dy_m,\preceq)$ is a distributive lattice called the \emph{Stanley lattice}~\cite{BeBICLRT}, \cite[Theorem~2.1]{FePLLP}.
Its bottom element is the sawtooth-shaped Dyck word $(\dw U\dw D)^m$, and its top element is the pyramid $\dw U^m\dw D^m$. The join ($\OR$) and meet ($\AND$) are computed pointwise on the height functions, and the two resulting height functions are called the upper and lower envelopes, respectively.

\smallskip
Given $w \in \Dy_m$, we define $\alpha_k(w)$ (for $k \in [m]$) as the number of $\dw U$'s appearing before the $(k+1)$-th occurrence of $\dw D$ in $w$. Note that $\alpha_k(w)\geq k+1$ for all $k \in [m]$ (the prefix of $w$ ending with the $(k+1)$-th occurrence of $\dw D$ contains $k+1$ occurrences of $\dw D$, hence at least as many of $\dw U$), $\alpha_{k+1}(w)-\alpha_k(w) \geq 0$ for all $k \in [m-1]$, and $\alpha_{m-1}(w)=m$.
The numbers $\alpha_k(w)$ are the partial sums of the \emph{dominating sequence} of $w$ of~\cite[pp.~2--3]{SaTTDPODP}, through which the dominance order is defined there.

\medskip
An \emph{$n$-staircase function} is a non-decreasing function $f: [n]\to [n+1]$ such that $f(x)\geq x$ for all $x \in [n]$. The set of all $n$-staircase functions is denoted by $\Stair_n$. 
Staircase functions can be ordered pointwise: $f\leq g$ iff $f(x)\leq g(x)$ for all $x \in [n]$.

\begin{theorem}
  \label{thm:rho}
  For $A \in \E_n$, $A \neq \top$, let $\rho_A: [n] \to [n+1]$ be the function defined by
  \[
    \rho_A(\ell) = \min\{\,y \mid [x\..y] \in A \text{ and } \ell \leq x\,\},
  \]
  with $\rho_A(\ell)=n$ if the set on the right-hand side is empty. Then:
  \begin{enumerate}
    \item $\rho_A$ is an $n$-staircase function;
    \item the map $A \mapsto \rho_A$ is an order isomorphism between $(\E_n\setminus\{\top\},\leq)$ and $(\Stair_n,\geq)$, where the order on $\Stair_n$ is the pointwise order.
  \end{enumerate}
\end{theorem}
\begin{proof}
\noindent 1. Suppose $\ell_1 \leq \ell_2$: then the set 
$\{\,y \mid [x\..y] \in A \text{ and } \ell_1 \leq x\,\}$ 
is a superset of $\{\,y \mid [x\..y] \in A \text{ and } \ell_2 \leq x\,\}$, so either the latter is empty (in which case $\rho_A(\ell_2)=n\geq \rho_A(\ell_1)$) or otherwise $\rho_A(\ell_1)\leq \rho_A(\ell_2)$ (the former being a minimum computed on a superset). Finally, $\rho_A(\ell)\geq \ell$ because either the set is empty (hence $\rho_A(\ell)=n\geq \ell$) or otherwise there exists $[x\..y]\in A$ such that $\ell\leq x$, so $y\geq x\geq \ell$ and $\rho_A(\ell)=y\geq \ell$.

\medskip
\noindent 2. Let us first prove \emph{monotonicity}: suppose $A \leq B$ and let $\ell \in [n]$. If the set on the right-hand side of the definition of $\rho_A(\ell)$ is empty, then $\rho_B(\ell)\leq n=\rho_A(\ell)$. Otherwise, let $[x\..y]\in A$ be such that $\ell\leq x$ and $y=\rho_A(\ell)$. By hypothesis, there exists $[x'\..y']\in B$ such that $x \leq x' \leq y' \leq y$, so $\ell \leq x \leq x'$ and $y' \leq y = \rho_A(\ell)$, hence $\rho_B(\ell)\leq y' \leq \rho_A(\ell)$.

To prove that the map \emph{reflects} the order (which also implies injectivity), let us show that $\rho_B \leq \rho_A$ implies $A \leq B$. Let $[x\..y]\in A$. By definition, $\rho_A(x)\leq y$,  and by hypothesis $\rho_B(x)\leq \rho_A(x)\leq y$, so there exists $[x'\..y']\in B$ such that $x \leq x' \leq y' \leq y$, hence $A \leq B$.

Finally, to prove \emph{surjectivity}, let $f\in\Stair_n$. We construct an antichain $A$ such that $\rho_A=f$. Let
\[
A=\sing{[x\..f(x)]\mid x\in[n],\ f(x)<n,\text{ and either $x=n-1$ or $f(x)<f(x+1)$}},
\]
that is, we keep, for each value of $f$, only the largest element of $[n]$ attaining it. $A$ is an antichain: if $[x\..f(x)],[y\..f(y)]\in A$ with $x<y$, then $f(x)<f(y)$ by construction, so neither interval contains the other.
We prove that $\rho_A=f$. Let $\ell\in[n]$. If $f(\ell)=n$, then every $x\geq\ell$ has $f(x)=n$ by monotonicity, so no interval of $A$ has left extreme at least $\ell$ and $\rho_A(\ell)=n$. Otherwise, let $t=\max\{\,x\mid f(x)=f(\ell)\,\}$: then $[t\..f(t)]\in A$, $t\geq\ell$ and $f(t)=f(\ell)$, so $\rho_A(\ell)\leq f(\ell)$; conversely, every $[x\..f(x)]\in A$ with $x\geq\ell$ has $f(x)\geq f(\ell)$ by monotonicity, so $\rho_A(\ell)=f(\ell)$.
\end{proof}

So there is an order isomorphism between $(\E_n\setminus\{\top\},\leq)$ and $(\Stair_n,\geq)$. On the other hand\footnote{Through the order-reversing bijection $a_1=1$, $a_i=n+1-f(n+1-i)$ ($2\leq i\leq n+1$), Theorem~\ref{thm:delta} is equivalent to the observation of~\cite[p.~320]{BBFPDLS} that increasing parking functions ordered componentwise form a lattice anti-isomorphic to the Stanley lattice.}:
\begin{theorem}
  \label{thm:delta}
  For $f\in\Stair_n$, let $a_i=f(i)+1$ for $i \in [n]$, let
  $a_{n}=n+1$, and define
  \[
    \delta(f)=\dw U^{a_0}\dw D\,\dw U^{a_1-a_0}\dw D\,\cdots\,
    \dw U^{a_n-a_{n-1}}\dw D.
  \]
  Then $\delta$ is an order isomorphism between $(\Stair_n,\leq)$ and $(\Dy_{n+1},\preceq)$. 
\end{theorem}
\begin{proof}
  First notice that all the superscripts appearing in the definition of $\delta(f)$ are natural numbers, because $a_{k+1}-a_k = f(k+1)-f(k) \geq 0$. Let us prove that $\delta(f)$ is a Dyck word of semilength $n+1$. The number of $\dw U$ appearing before the $(k+1)$-th occurrence of $\dw D$ is
  \[
    a_0 + (a_1-a_0) + \cdots + (a_k-a_{k-1}) = a_k = f(k)+1 \geq k+1,
  \]
  because $f(k)\geq k$: hence no prefix contains more $\dw D$'s than $\dw U$'s. Moreover there are $n+1$ occurrences of $\dw D$, and $a_n= n+1$ occurrences of $\dw U$.

  From $\delta(f)$ we can fully recover $f$, so $\delta$ is injective. To prove surjectivity, consider a Dyck word $w$ of semilength $n+1$. Define $f_w(k)=\alpha_k(w)-1$ for $k \in [n]$. Then, given the known properties of $\alpha_k(w)$, we have $f_w\in\Stair_n$ and $\delta(f_w)=w$.

  Finally, let us prove that $\delta$ preserves and reflects the order (a bijective order-preserving map between posets need not be an isomorphism, so both directions are required). Writing $d_w(t)$ for the number of $\dw D$'s in the prefix of length $t$ of $w$, we have $h_w(t)=t-2d_w(t)$, so for $w,v \in \Dy_{n+1}$
  \begin{multline*}
    w \preceq v \text{ iff } 
    d_w(t)\geq d_v(t)\ \text{for all $t$}\\
    \text{iff each occurrence of $\dw D$ appears in $w$ no later than in $v$,}
  \end{multline*}
  and the position of the $(k+1)$-th occurrence of $\dw D$ in $w$ is $\alpha_k(w)+k+1$: hence $w \preceq v$ iff $\alpha_k(w)\leq\alpha_k(v)$ for all $k \in [n+1]$. Since $\alpha_k(\delta(f))=f(k)+1$ for $k \in [n]$ and $\alpha_n(\delta(f))=n+1$, we conclude that $\delta(f)\preceq\delta(g)$ iff $f\leq g$.
\end{proof}

\begin{figure}
	\centering
	\includegraphics[scale=.3]{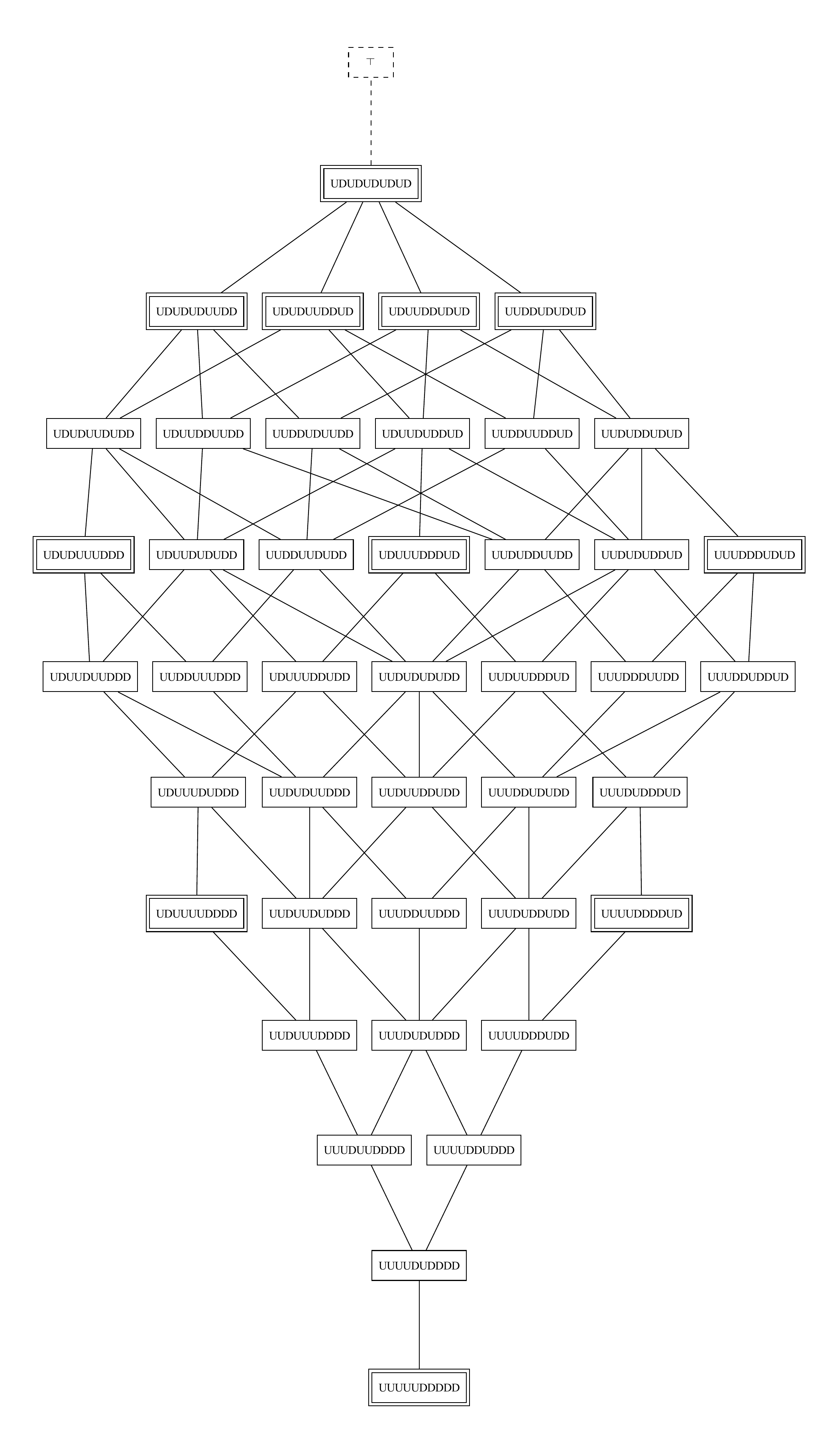}
	\caption{\label{fig:e4-dyck}The lattice $\E_4 \setminus \{\,\top\,\}$ of Figure~\ref{fig:e4}, but every antichain $A$ is mapped to the Dyck word $\delta(\rho_A)$ (see Theorem~\ref{thm:rho} and Theorem~\ref{thm:delta}). Elements with double frame are $\AND$-irreducible.}
\end{figure}

Composing the results obtained so far, we have the following chain of order isomorphisms:
\[
\begin{tikzcd}[row sep=small]
  (\E_n\setminus\{\top\},\leq) \arrow[r] & (\Stair_n,\geq) \arrow[r] & (\Dy_{n+1},\succeq) \\[-4pt]
  A \arrow[r,mapsto,shorten <=9mm,shorten >=4mm] & \rho_A \arrow[r,mapsto,shorten <=4mm,shorten >=4mm] & \delta(\rho_A)
\end{tikzcd}
\]
(see also Section~\ref{sec:related} for the relation with Birkhoff's theorem and~\cite[Proposition~3.7]{FeMLPRTV}).
Finally, it is useful to define:
\begin{eqnarray*}
  \omega: (\E_n\setminus\{\top\},\leq) & \longrightarrow & (\Dy_{n+1},\succeq) \\  
  A & \longmapsto & \omega(A)=\delta(\rho_A).
\end{eqnarray*}
Thus $\omega$ is an \emph{anti-isomorphism} between $\E_n\setminus\{\top\}$ and the Stanley lattice $(\Dy_{n+1},\preceq)$, that is, an order isomorphism onto its dual; we use this terminology throughout the paper.

\begin{figure}
  \centering
  \resizebox{0.79\textwidth}{!}{\input{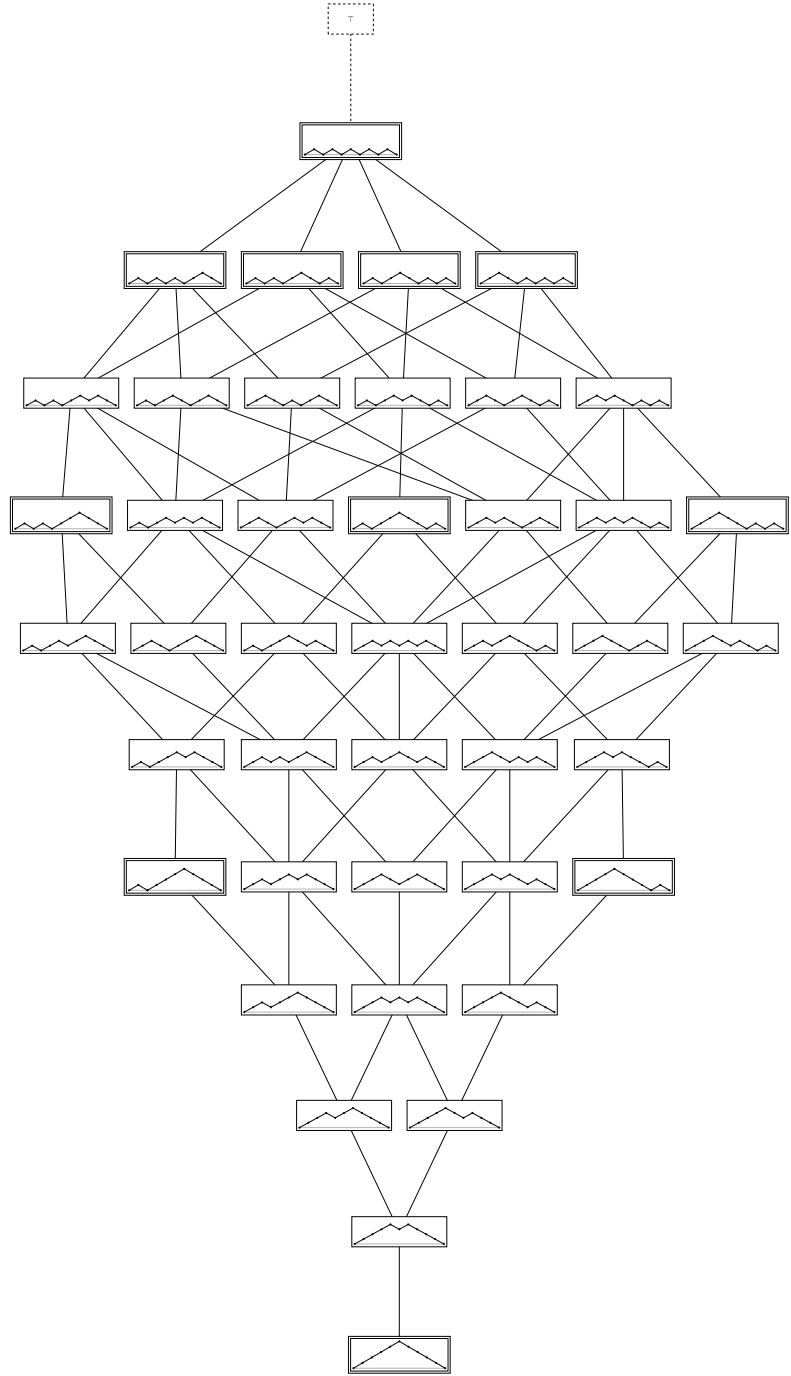}}
  \caption{\label{fig:e4-dyck-h}The lattice $\E_4 \setminus \{\,\top\,\}$ of Figure~\ref{fig:e4-dyck}, but now every Dyck word $w$ is represented by its height function $h_w$: observe that order represents reverse dominance (lower elements in the poset dominate upper elements). Elements with double frame are $\AND$-irreducible.}
\end{figure}

\subsection{Area and rank}

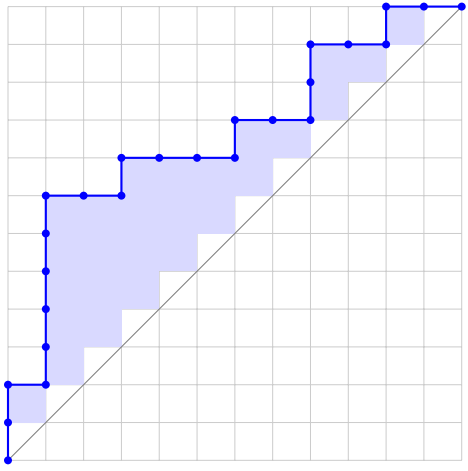
\begin{figure}
  \centering
  \begin{tikzpicture}[scale=.5]
  \draw[very thin, gray!40] (0,0) grid (12,12);
  \fill[blue!15] (0,1) rectangle (1,2);
  \fill[blue!15] (1,2) rectangle (2,3);
  \fill[blue!15] (1,3) rectangle (2,4);
  \fill[blue!15] (1,4) rectangle (2,5);
  \fill[blue!15] (1,5) rectangle (2,6);
  \fill[blue!15] (1,6) rectangle (2,7);
  \fill[blue!15] (2,3) rectangle (3,4);
  \fill[blue!15] (2,4) rectangle (3,5);
  \fill[blue!15] (2,5) rectangle (3,6);
  \fill[blue!15] (2,6) rectangle (3,7);
  \fill[blue!15] (3,4) rectangle (4,5);
  \fill[blue!15] (3,5) rectangle (4,6);
  \fill[blue!15] (3,6) rectangle (4,7);
  \fill[blue!15] (3,7) rectangle (4,8);
  \fill[blue!15] (4,5) rectangle (5,6);
  \fill[blue!15] (4,6) rectangle (5,7);
  \fill[blue!15] (4,7) rectangle (5,8);
  \fill[blue!15] (5,6) rectangle (6,7);
  \fill[blue!15] (5,7) rectangle (6,8);
  \fill[blue!15] (6,7) rectangle (7,8);
  \fill[blue!15] (6,8) rectangle (7,9);
  \fill[blue!15] (7,8) rectangle (8,9);
  \fill[blue!15] (8,9) rectangle (9,10);
  \fill[blue!15] (8,10) rectangle (9,11);
  \fill[blue!15] (9,10) rectangle (10,11);
  \fill[blue!15] (10,11) rectangle (11,12);
  \draw[gray] (0,0) -- (12,12);
  \draw[thick, blue] (0,0) -- (0,1) -- (0,2) -- (1,2) -- (1,3) -- (1,4) -- (1,5) -- (1,6) -- (1,7) -- (2,7) -- (3,7) -- (3,8) -- (4,8) -- (5,8) -- (6,8) -- (6,9) -- (7,9) -- (8,9) -- (8,10) -- (8,11) -- (9,11) -- (10,11) -- (10,12) -- (11,12) -- (12,12);
  \fill[blue] (0,0) circle (3pt);
  \fill[blue] (0,1) circle (3pt);
  \fill[blue] (0,2) circle (3pt);
  \fill[blue] (1,2) circle (3pt);
  \fill[blue] (1,3) circle (3pt);
  \fill[blue] (1,4) circle (3pt);
  \fill[blue] (1,5) circle (3pt);
  \fill[blue] (1,6) circle (3pt);
  \fill[blue] (1,7) circle (3pt);
  \fill[blue] (2,7) circle (3pt);
  \fill[blue] (3,7) circle (3pt);
  \fill[blue] (3,8) circle (3pt);
  \fill[blue] (4,8) circle (3pt);
  \fill[blue] (5,8) circle (3pt);
  \fill[blue] (6,8) circle (3pt);
  \fill[blue] (6,9) circle (3pt);
  \fill[blue] (7,9) circle (3pt);
  \fill[blue] (8,9) circle (3pt);
  \fill[blue] (8,10) circle (3pt);
  \fill[blue] (8,11) circle (3pt);
  \fill[blue] (9,11) circle (3pt);
  \fill[blue] (10,11) circle (3pt);
  \fill[blue] (10,12) circle (3pt);
  \fill[blue] (11,12) circle (3pt);
  \fill[blue] (12,12) circle (3pt);
\end{tikzpicture}
  \caption{\label{fig:dyckarea}The Dyck word $w=\dw{UUDUUUUUDDUDDDUDDUUDDUDD}$ of Figure~\ref{fig:dyckword}, drawn as a lattice path from $(0,0)$ to $(12,12)$ with $\dw U$ as a north step and $\dw D$ as an east step. The shaded cells are the complete lattice squares between the path and the diagonal: $\area(w)=26$.}
\end{figure}

The \emph{area} of a Dyck word $w \in \Dy_m$ is the number of complete lattice squares between the path and the diagonal, when $w$ is drawn as a lattice path from $(0,0)$ to $(m,m)$ with north and east steps staying weakly above the diagonal (see Figure~\ref{fig:dyckarea}); this is the standard area statistic on Dyck paths~\cite{HagQTC}. It can be computed as\footnote{An equivalent formula is~\cite[equation (1)]{SaTTDPODP}; moreover, the observation that the area of a Dyck word is its rank in the Stanley lattice is~\cite[Proposition~3.1]{FePLLP}.}
\[
  \area(w)=\sum_{i \in [m]} \alpha_i(w) - \binom{m+1}{2}.
\]

The maximum area of a Dyck word of semilength $m$ is $\binom{m}{2}$~\cite[p.~5]{SaTTDPODP}, and it is attained by the word $w=\dw U^m\dw D^m$: in fact, for this word $\alpha_k(w)=m$ for all $k \in [m]$, so the area is $m^2-\binom{m+1}{2}=\binom{m}{2}$. 

The minimum area is $0$, and it is attained by the word $v=(\dw U\dw D)^m$: in fact, for this word $\alpha_k(v)=k+1$ for all $k \in [m]$, so the area is $\binom{m+1}{2} - \binom{m+1}{2} = 0$.

\medskip
The \emph{rank} of an element $A \in \E_n$, denoted by $\rank(A)$, is the length of the longest chain from $\bot$ to $A$. In fact, all maximal chains have the same length, which is the rank of $A$. Furthermore, in~\cite{BoVLAFI} it is shown that the rank of $A$ can be computed as follows:

\begin{theorem}
  \label{thm:rank}
  Let $[\ell\..r]\in \I_n$; then,
\[
\rank(\{\?[\ell\..r]\?\}) = (1+\ell)(n-r).
\]
More generally, let $k\geq0$ and $A=\{\? [\ell_0\..r_0],\dots,[\ell_k\..r_k] \?\}\in\E_n$ with
$\ell_0<\ell_1<\dots<\ell_k$ (so that $A\neq\emptyset$ and $A\neq\top$). Then,
\[
	\rank(A)=(1+\ell_0)(n-r_0)+\sum_{i=1}^k (\ell_i-\ell_{i-1})(n-r_i),
\]
of which the previous formula is the case $k=0$ (the sum being then empty). Finally, $\rank(\emptyset)=0$.
\end{theorem}

There is a simple relationship between the area of a Dyck word and the rank of the corresponding antichain:
\begin{theorem}
  \label{thm:area-rank}
  Let $A \in \E_n$, $A \neq \top$. Then,
  \[
    \area(\omega(A)) = \binom{n+1}{2} - \rank(A).
  \]
\end{theorem}
\begin{proof}
  By definition of $\delta$ (see the proof of Theorem~\ref{thm:delta}), $\alpha_k(\omega(A))=\rho_A(k)+1$ for $k \in [n]$, and $\alpha_n(\omega(A))=n+1$. Hence
  \begin{multline*}
    \area(\omega(A)) = \sum_{k \in [n+1]} \alpha_k(\omega(A)) - \binom{n+2}{2} =\\
    = \sum_{k \in [n]} \rho_A(k) + 2n+1 - \binom{n+2}{2}.
  \end{multline*}
  We now claim that
  \[
    \rank(A) = \sum_{x \in [n]} \bigl(n - \rho_A(x)\bigr).
  \]
  If $A=\emptyset$ (the bottom of $\E_n$), both sides are $0$, as $\rho_A$ is constantly $n$. Otherwise, write $A=\{\?[\ell_0\..r_0],\dots,[\ell_k\..r_k]\?\}$ with $\ell_0<\dots<\ell_k$; since $A$ is an antichain, also $r_0<\dots<r_k$. Then $\rho_A(x)=r_0$ for $x\leq\ell_0$, $\rho_A(x)=r_i$ for $\ell_{i-1}<x\leq\ell_i$, and $\rho_A(x)=n$ for $x>\ell_k$, so
  \[
    \sum_{x \in [n]} \bigl(n-\rho_A(x)\bigr)
    = (1+\ell_0)(n-r_0)+\sum_{i=1}^k (\ell_i-\ell_{i-1})(n-r_i) = \rank(A)
  \]
  by Theorem~\ref{thm:rank}. Adding the two identities,
  \begin{multline*}
    \area(\omega(A))+\rank(A) = n^2+2n+1- \binom{n+2}{2}=\\
    =(n+1)^2- \binom{n+2}{2}=\binom{n+1}{2}.
  \end{multline*}
\end{proof}
The same result can be derived from the fact that $\omega$ is an anti-isomorphism between graded lattices, and the fact that the rank of a Dyck word in the Stanley lattice is its area~\cite[Proposition~3.1]{FePLLP}.

From the above result, we have the following corollary:
\begin{corollary}
  \label{cor:rank-area}
  Let $A,A' \in \E_n \setminus \{\,\top\,\}$.
  Then $\rank(A) = \rank(A')$ iff $\area(\omega(A)) = \area(\omega(A'))$. In other words, the rank of an antichain is determined by the area of the corresponding Dyck word.
\end{corollary}

\paragraph{Area and height function.} The definition we used for the area of a Dyck word is standard in the Dyck path literature, but it is not the only one: the area of a Dyck word $w \in \Dy_m$ is sometimes defined as the sum of its heights $\sum_{t \in [2m+1]} h_w(t)$. The two notions are related as follows~\cite[Proposition~3.1]{FePLLP}:
\begin{equation}
  \label{eqn:area-height}
  \area(w) = \frac{1}{2}\left(\sum_{t\in [2m+1]} h_w(t) - m\right).
\end{equation}

\subsection{Carlitz--Riordan $q$-Catalan numbers and rank generating function}

Recall that the Carlitz--Riordan $q$-Catalan number~\cite{CaRTELPNTQG} $\overline C_n(q)$ is a polynomial in $q$ defined by
\[
  \overline C_n(q) = \sum_{w \in \Dy_n} q^{\area(w)}.
\]
Areas of Dyck words of semilength $n$ range from $0$ to $\binom{n}{2}$, so $\overline C_n(q)$ can be rewritten as
\[
  \overline C_n(q) = \sum_{t=0}^{\binom{n}{2}} a_n(t) q^t,
\]
where $a_n(t)$ is the number of Dyck words of semilength $n$ with area $t$. (These polynomials should not be confused with MacMahon's $q$-Catalan numbers $\frac1{[n+1]_q}\binom{2n}{n}_q$, a different $q$-analogue of the Catalan numbers; see~\cite{FuHQCN} for both.) By Corollary~\ref{cor:rank-area}, the level sets of $\E_n\setminus\{\top\}$ (the sets of elements of the same rank) are in bijection with the sets of Dyck words of semilength $n+1$ of the same area, so their sizes are the coefficients of $\overline C_{n+1}(q)$. In particular, the size of the largest \emph{level set} of $\E_n$ is the maximum coefficient of the Carlitz--Riordan $q$-Catalan polynomial $\overline C_{n+1}(q)$ (the top element forms a level set of size $1$ by itself, which does not affect the maximum). 

\subsubsection{Level set sizes and integer partitions}\label{sec:partitions}

Let $\lambda_n$ be the maximum coefficient of $\overline C_{n+1}(q)$, that is, the size of the largest level set of $\E_n$.
Since level sets are antichains, we have that
\[
  \lambda_n \leq \width(\E_n); 
\]
whether equality holds (that is, whether $\E_n$ enjoys the Sperner property) is an open problem. 
It was verified in~\cite{BoVLAFI} for $n\leq11$, and it holds for $n\leq20$: the widths of the Stanley lattices form sequence A274291 of the OEIS~\cite{OEISA274291}, which is indexed by the semilength, so that $\width(\E_n)$ is its term of index $n+1$, and the terms listed in the entry, up to index $21$, coincide with the values of $\lambda_n$ for $n\leq20$ of Table~\ref{tab:lambda}.
If the Sperner property holds for all $n$, the open problem of determining the width of $\E_n$ is equivalent to that of computing $\lambda_n$. 
The table up to semilength $300$ computed by T.~M\"utze and linked from the entry is explicitly computed under the assumption of the Sperner property, that is, it lists the values of $\lambda_n$. The values of $\lambda_n$ for $n\leq20$ also appear, as the largest coefficients of the $q$-Catalan polynomials, in Table~I of Winston and Kleitman~\cite{WiKANTSS}.

Recall that, by Theorem~\ref{thm:rho},  every $A \in \E_n\setminus\{\top\}$ is mapped to an $n$-staircase function $\rho_A$ (a non-decreasing function such that $\rho_A(x)\geq x$), and any staircase function $f \in \Stair_n$ can be identified with the set of cells\footnote{We are identifying the pair $(x,y)\in [n]^2$ with the cell $[x,x+1] \times [y,y+1]$ in the $n \times n$ grid.}
\[
  S_f = \{\? (x, n-y-1) \in [n]^2 \mid y \geq f(x) \?\},
\]
a lower set of the poset $([n]^2,\leq)$, contained in the triangle 
\[
  T_n = \{\? (x,y) \in [n]^2 \mid x+y \leq n-1 \?\}.
\]
Since $\rank(A)=\sum_{x \in [n]} (n-\rho_A(x))$, the rank of $A$ is the number of cells in $S_{\rho_A}$. The column heights $\mu_{x+1}=n-\rho_A(x)$ are non-increasing in $x$, so $\mu=(\mu_1,\dots,\mu_n)$ is an integer partition of $\rank(A)$, and the Ferrers diagram~\cite{StaEC1} of $\mu$ (drawn with parts as column heights, justified at the bottom) is contained in $T_n$, which is the Ferrers diagram of the partition $\mathsf{st}_n=(n,n-1,\dots,1)$. Conversely, every integer partition $\mu$ whose Ferrers diagram is contained in $T_n$ arises in this way.
So
\begin{equation}
  \label{eqn:rank-partition}
  \sum_{A \in \E_n\setminus\{\top\}} q^{\rank(A)} = \sum_{\mu \subseteq \mathsf{st}_n} q^{|\mu|}.
\end{equation}
By Theorem~\ref{thm:area-rank}, for each $A$ we have
\[
  q^{\rank(A)} = q^{\binom{n+1}{2}-\area(\omega(A))} = q^{\binom{n+1}{2}}\left(\frac1q\right)^{\area(\omega(A))}. 
\]
Summing over all $A \in \E_n\setminus\{\top\}$, we obtain
\begin{multline*}
  \sum_{A \in \E_n\setminus\{\top\}} q^{\rank(A)} = q^{\binom{n+1}{2}} \sum_{A \in \E_n\setminus\{\top\}} \left(\frac1q\right)^{\area(\omega(A))} =\\
  = q^{\binom{n+1}{2}} \sum_{w \in \Dy_{n+1}} \left(\frac1q\right)^{\area(w)} = q^{\binom{n+1}{2}} \overline C_{n+1}\left(\frac1q\right). 
\end{multline*}
Now using~\eqref{eqn:rank-partition}, we obtain:
\[
  \sum_{\mu \subseteq \mathsf{st}_n} q^{|\mu|} = q^{\binom{n+1}{2}} \overline C_{n+1}\left(\frac1q\right) = \overline C_{n+1}^*(q),
\]
where $\overline C_{n+1}^*(q)$ is the \emph{reverse} of the polynomial $\overline C_{n+1}(q)$: the latter equality follows because
 $\overline C_{n+1}(q)$ is a polynomial of degree $\binom{n+1}{2}$. Note that $\overline C_{n+1}^*(q)$ has the same coefficients as $\overline C_{n+1}(q)$, but in reverse order.

Hence, for $0 \leq t \leq \binom{n+1}{2}$, the coefficient $a_{n+1}(t)$ of $q^t$ in $\overline C_{n+1}(q)$ can be identified with the number of integer partitions $\mu$ whose Ferrers diagram is contained in $T_n$ and such that $|\mu|=\binom{n+1}{2} - t$:
\[
  a_{n+1}(t) = \#\Bigl\{\? \mu \subseteq \mathsf{st}_n \mid |\mu|=\binom{n+1}{2} - t \?\Bigr\}.
\]
(The restriction on $t$ matters: for $t>\binom{n+1}{2}$ the left-hand side is $0$, whereas the right-hand side, read with truncated subtraction, would count the empty partition.) This identity is well known~\cite[Proposition~3.9]{FeMLPRTV}, \cite[p.~3]{FerUDP}, \cite[p.~52]{StaUYL}; the derivation above shows how it arises from $\E_n$.

\smallskip
This observation connects the coefficients of the Carlitz--Riordan $q$-Catalan numbers with integer partitions contained in $\mathsf{st}_n$, which suggests a recursive approach. For each positive integer $m$ and each $j \in [m]$, let $P_{m,j}(q)$ be the generating function of integer partitions contained in $\mathsf{st}_{m-1}$ and whose first part is at most $j$. Removing the first part of such a partition (the leftmost column of its Ferrers diagram) produces a partition contained in $\mathsf{st}_{m-2}$ whose first part is at most the removed one. Hence, we have
\begin{eqnarray*}
  P_{m,0}(q) & = & 1\quad(m\geq1),\\
  P_{m,j}(q) & = & P_{m,j-1}(q)+ q^j P_{m-1,\min(j,m-2)}(q).
\end{eqnarray*}
Finally, we have:
\[
  \overline C_{n+1}^*(q)=\sum_{\mu \subseteq \mathsf{st}_n} q^{|\mu|} = P_{n+1,n}(q).
\]
In particular, $\lambda_n$ is the maximum coefficient of $P_{n+1,n}(q)$, which can be computed recursively. The polynomials $\overline C_{n+1}(q)$ can also be computed by the classical convolution $\overline C_{k+1}(q)=\sum_{\nu=0}^{k}q^\nu\,\overline C_\nu(q)\,\overline C_{k-\nu}(q)$~\cite[Proposition~3.2]{SaTTDPODP}, so the only merit of the recursion above is its cost: it computes $\overline C^*_{n+1}(q)$ with $O(n^4)$ integer additions, whereas the convolution requires $O(n^6)$ integer multiplications if polynomials are multiplied naively.

\subsection{An asymptotic estimate of $\lambda_n$}\label{sec:asymptotics}

The recursive formula for $P_{m,j}(q)$ can be used to compute $\lambda_n$ exactly, but it is not a closed formula, and it says nothing about the growth of $\lambda_n$.
We can, however, provide an asymptotic estimate of $\lambda_n$ using the connection between Dyck words and random walks. Let us recall some terminology.

Let $X_1, X_2, \dots$ be a sequence of independent and identically distributed random variables with $\Pr(X_i=-1)=\Pr(X_i=1)=1/2$, 
and let $S_k = X_1 + \dots + X_k$. The sequence $S_0, S_1, \dots$ is called a \emph{simple symmetric random walk}: its values are integers, and the walk starts at $0$ and moves up or down by $1$ with equal probability at each step.
Following~\cite{DKWLAARWE}, let $\tau=\min\{\,k\geq1\mid S_k\leq0\,\}$ be the first time the walk is nonpositive. Throughout this subsection, $m=n+1$ denotes the semilength of the Dyck words associated with $\E_n$. If we condition the random walk to the event $\tau=2m+2$, the walk stays strictly positive for $2m+1$ steps and returns to $0$ at step $2m+2$; hence $S_1=S_{2m+1}=1$, and the sequence $S_1-1,\dots,S_{2m+1}-1$ is the height function of a Dyck word $w$ of semilength $m$. Since every such path has probability $2^{-(2m+2)}$, the Dyck word $w$ is extracted uniformly at random. More precisely, $S_k=h_w(k-1)+1$ for $1\leq k\leq 2m+1$; if we further define $A_{2m+1}=S_1+\dots+S_{2m+1}$, using~\eqref{eqn:area-height} we have
\[
  A_{2m+1}=\sum_{t\in[2m+1]}h_w(t)+2m+1=2\,\area(w)+3m+1.
\]
Since we already know that $a_m(t)$ is the number of Dyck words of semilength $m$ with area $t$, we have that extracting a Dyck word of semilength $m$ uniformly at random, the probability that it has area $t$ is
\[
\Pr\left[\area(w)=t\right] = a_m(t) / C_m
\]
(where $C_m$ is the $m$-th Catalan number, the total number of Dyck words of semilength $m$). 
By the observations above, this quantity can be rewritten as
\begin{equation}
  \label{eqn:area-prob}
  \Pr\left[A_{2m+1}=2t+3m+1 \mid \tau=2m+2\right] = a_m(t) / C_m.
\end{equation}
Now, to estimate the maximum value on the left-hand side, we use the following local limit theorem, stated here in the notation of its authors, except that we write $\varrho$ and $N$ for their $\rho$ and $n$ (in particular, $\varrho$ is a parameter of the step distribution, unrelated to our $\rho_A$).
\begin{theorem}[{Denisov, Kolb and Wachtel~\cite[Theorem~1.1, first assertion]{DKWLAARWE}}]
  \label{thm:dkw}
  Let $S_k=X_1+\dots+X_k$, where the $X_i$ are independent, identically distributed and integer-valued, with $\mathbf E X_i=0$ and $\mathbf E X_i^2=\sigma^2\in(0,\infty)$, and assume that $X_i$ is $(d,\varrho)$-lattice, that is, $d$ is the largest integer such that $\Pr(X_i\in\{\?\varrho+dk\mid k\in\mathbf Z\?\})=1$, with $\varrho\in[0,d)$. Let $\tau=\min\{\,k\geq1\mid S_k\leq0\,\}$ and $A_N=S_1+\dots+S_N$. Then
  \[
    \sup_{a\in N(N+1)\varrho/2+d\mathbf Z}\left|N^{3/2}\Pr\left[A_N=a\mid\tau=N+1\right]-\frac d\sigma\,w_{\mathrm{ex}}\!\left(\frac{a}{\sigma N^{3/2}}\right)\right|\longrightarrow0\qquad\text{as $N\to\infty$,}
  \]
  where $w_{\mathrm{ex}}$ is the density of the random variable $\int_0^1e(t)\,dt$, with $e$ a standard Brownian excursion.
\end{theorem}
In our setting the steps are $\pm1$, so $\sigma^2=1$, $d=2$ and $\varrho=1$ (the steps lie in $1+2\mathbf Z$), and $d/\sigma=2$. The excursion we condition on is $S_1,\dots,S_{2m+1}$ ($\tau=2m+2$), so we apply the theorem with $N=2m+1$. The condition on $a$ in the supremum, $a\in(2m+1)(m+1)+2\mathbf Z$, is automatically satisfied by $a=2t+3m+1$, as $(2m+1)(m+1)=2m^2+3m+1$, so the statement of the theorem reduces to the following:
\begin{multline*}
  \sup_t \left| (2m+1)^{3/2} \cdot \Pr\left[A_{2m+1}=2t+3m+1 \mid \tau=2m+2\right] \right. \\\left.- 2 \psi\left(\frac{2t+3m+1}{(2m+1)^{3/2}}\right)\right| \to 0 \quad \text{as $m \to \infty$},
\end{multline*}
where $\psi=w_{\mathrm{ex}}$ is the density of the Airy distribution~\cite{TakBEVA} (the law of the area under a Brownian excursion), a bounded continuous function.
Hence (recalling that $\lambda_n$ is the maximum coefficient of $\overline C_{n+1}(q)=\overline C_m(q)$) we have
\begin{multline*}
  \lambda_n = \max_t a_m(t) = \text{ (using~\eqref{eqn:area-prob})}\\ = C_m \cdot \max_t \Pr\left[A_{2m+1}=2t+3m+1 \mid \tau=2m+2\right] \sim\\
  \sim \frac{2 C_m}{(2m+1)^{3/2}} \sup_x \psi(x)\sim \frac{C_m}{\sqrt2\, m^{3/2}} \sup_x \psi(x),
\end{multline*}
where the asymptotic equivalence holds because the error term in the theorem is uniform in $t$ and $\psi$ is continuous: for $0\leq t\leq\binom m2$ the points $(2t+3m+1)/(2m+1)^{3/2}$ form a grid whose mesh tends to $0$, whose first point tends to $0$ and whose last point tends to $\infty$, so the maximum of $\psi$ over the grid tends to $\sup_x \psi(x)$, which is a constant (it does not depend on $m$).
Using $C_m \sim \frac{4^m}{\sqrt{\pi}\, m^{3/2}}$, we obtain
\[
  \lambda_n\sim \frac{4^m}{\sqrt{2\pi}\, m^3} \sup_x \psi(x) = \kappa\,\frac{4^{n+1}}{(n+1)^3},
\]
where $\kappa=\sup_x \psi(x) / \sqrt{2\pi}$ is a constant, whose value is approximately $1.080744$. Winston and Kleitman~\cite{WiKANTSS}, who met these coefficients while counting tournament score sequences, conjectured that the largest coefficient of $\overline C_m(q)$ is of the order of $C_m$ times the maximum of a normal density with the same variance as the distribution of the area, that is, $O(4^m/m^3)$~\cite[Section~9]{WiKANTSS}; this upper bound was proved by Kim and Pittel~\cite{KiPCKWC}, whereas a matching lower bound follows from Chebyshev's inequality, as already observed by Winston and Kleitman~\cite[(9.5)]{WiKANTSS}. Tak\'acs~\cite[(39)]{TakBEVA} showed that, if the coefficients of $\overline C_m(q)$ are unimodal and the limit theorem for the area can be replaced by a local limit theorem, then the largest coefficient of $\overline C_m(q)$ is asymptotic to $\kappa\,4^m/m^3$, with $\kappa=\psi(0.563\ldots)/\sqrt{2\pi}=1.080744143\ldots$, where $0.563\ldots$ is the mode of the Airy density $\psi$ (Tak\'acs denotes the density by $w$ and the distribution function by $W$); he also noted that the ratio between the largest coefficient and the normal estimate of Winston and Kleitman then tends to $1.051271514\ldots$~\cite[(41)]{TakBEVA}, a value consistent with the ratios tabulated by Winston and Kleitman, which they observed to be ``a little greater than one''~\cite[Section~9]{WiKANTSS}. (Tak\'acs attributes to Winston and Kleitman the stronger conjecture that this limit is $1$~\cite[after (41)]{TakBEVA}; their Conjecture~1 only asserts that the ratio between the largest coefficient and the normal estimate is bounded, and the limit of that ratio is $1.0513\ldots$ by Tak\'acs's computation, which the theorem proved here makes unconditional, hence not $1$.) The local limit theorem of~\cite{DKWLAARWE} makes Tak\'acs's result unconditional; as shown above, unimodality is not needed, because the error term of the local limit theorem is uniform.

Table~\ref{tab:lambda} compares $\lambda_n$ with the asymptotic estimate for
$n\leq 20$, and Table~\ref{tab:lambda-large} does the same for larger values of
$n$. The convergence of the ratio is slow. The local limit theorem provides
only an $o(1)$ error term, hence no rate; empirically, however, the relative
error of the estimate behaves like $c/(n+1)$: the product
$(n+1)\bigl(1-\lambda_n/(\kappa\,4^{n+1}/(n+1)^3)\bigr)$ is about $1.2$ at
$n=20$ and grows slowly to about $1.28$ at $n=299$.

\begin{table}
  \centering
  \begin{tabular}{rrrr}
    \toprule
    $n$ & $\lambda_n$ & $\kappa\,4^{n+1}/(n+1)^3$ & ratio \\
    \midrule
    1 & \num{1} & \num{2} & 0.4626 \\
    2 & \num{2} & \num{3} & 0.7807 \\
    3 & \num{3} & \num{4} & 0.6940 \\
    4 & \num{7} & \num{9} & 0.7907 \\
    5 & \num{17} & \num{20} & 0.8295 \\
    6 & \num{44} & \num{52} & 0.8523 \\
    7 & \num{118} & \num{138} & 0.8530 \\
    8 & \num{338} & \num{389} & 0.8697 \\
    9 & \num{1003} & \num{1133} & 0.8851 \\
    10 & \num{3039} & \num{3406} & 0.8923 \\
    11 & \num{9466} & \num{10493} & 0.9021 \\
    12 & \num{30009} & \num{33012} & 0.9090 \\
    13 & \num{96757} & \num{105725} & 0.9152 \\
    14 & \num{316429} & \num{343834} & 0.9203 \\
    15 & \num{1047683} & \num{1133242} & 0.9245 \\
    16 & \num{3511473} & \num{3779166} & 0.9292 \\
    17 & \num{11876457} & \num{12734597} & 0.9326 \\
    18 & \num{40537388} & \num{43311370} & 0.9360 \\
    19 & \num{139490014} & \num{148536344} & 0.9391 \\
    20 & \num{483393651} & \num{513245115} & 0.9418 \\
    \bottomrule
  \end{tabular}
  \caption{\label{tab:lambda}The size $\lambda_n$ of the largest level set of $\E_n$ (the maximum coefficient of $\overline C_{n+1}(q)$), the asymptotic estimate $\kappa\,4^{n+1}/(n+1)^3$ (rounded to the nearest integer), and their ratio.}
\end{table}

\begin{table}
  \centering
  \begin{tabular}{rrrr}
    \toprule
    $n$ & $\lambda_n$ & $\kappa\,4^{n+1}/(n+1)^3$ & ratio \\
    \midrule
    25 & \num{2.6379e+11} & \num{2.7693e+11} & 0.9526 \\
    30 & \num{1.6061e+14} & \num{1.6730e+14} & 0.9600 \\
    40 & \num{7.3511e+19} & \num{7.5828e+19} & 0.9694 \\
    50 & \num{4.0292e+25} & \num{4.1312e+25} & 0.9753 \\
    60 & \num{2.4792e+31} & \num{2.5316e+31} & 0.9793 \\
    80 & \num{1.1702e+43} & \num{1.1889e+43} & 0.9843 \\
    100 & \num{6.6575e+54} & \num{6.7425e+54} & 0.9874 \\
    150 & \num{2.5361e+84} & \num{2.5577e+84} & 0.9915 \\
    200 & \num{1.3659e+114} & \num{1.3747e+114} & 0.9936 \\
    250 & \num{8.9030e+143} & \num{8.9487e+143} & 0.9949 \\
    299 & \num{1.6539e+173} & \num{1.6609e+173} & 0.9957 \\
    \bottomrule
  \end{tabular}
  \caption{\label{tab:lambda-large}The same comparison as in Table~\ref{tab:lambda} for larger values of $n$ (values rounded to five significant digits); the values of $\lambda_n$ are taken from the table of sequence A274291 in the OEIS~\cite{OEISA274291}, computed by T.~M\"utze.}
\end{table}

\section{Extensions: antichain lattice of intervals of a poset}\label{sec:ext}

We now move to the more general case of the antichain lattice of an \emph{arbitrary} finite poset $(P,\sqleq)$. The parameter that will govern the picture is the \emph{sweep dimension}, introduced below; by Proposition~\ref{prop:dimsweep} it is the Dushnik--Miller dimension of the lattice of order-convex subsets of $P$ (the intervals of $P$ in our sense), which connects the results of this section to dimension theory and to the theory of convex geometries (Section~\ref{sec:related}).

\subsection{Dimension and sweep dimension}

Every linear order that extends a partial order $\sqleq$ is called a \emph{linear extension} of $\sqleq$. 
The \emph{dimension} of a finite poset $(P,\sqleq)$ is the minimum number of linear extensions whose intersection is $\sqleq$, and it is denoted by $\dim(P,\sqleq)$~\cite[\S2.2]{DuMPOS} (see also~\cite{TroCPOSET}). We shall use the following reformulation, which is a variant of the one given by Dushnik and Miller~\cite[\S2.4]{DuMPOS} (where each coordinate is required to be injective and the order is strict); we include a proof for completeness.

\begin{lemma}
  \label{lem:dimembed}
  A finite poset $(P,\sqleq)$ has dimension at most $k$ if and only if there is an order embedding of $(P,\sqleq)$ into $\R^k$ with the componentwise order.
\end{lemma}
\begin{proof}
  If $\sqleq$ is the intersection of the linear extensions $L_1$,~\dots,~$L_k$ (not necessarily distinct: we repeat one of them if the dimension is smaller than $k$), mapping $x$ to the vector of its positions in $L_1$,~\dots,~$L_k$ is an order embedding. Conversely, let $g:P\to\R^k$ be an order embedding, and fix a linear extension $L$ of $\sqleq$. For each $i$, let $L_i$ order $P$ by the $i$-th coordinate of $g$, breaking ties using $L$. If $x\sqleq y$, then $g(x)_i\leq g(y)_i$, and in case of equality $x$ precedes $y$ in $L$, so $L_i$ is a linear extension of $\sqleq$. If $x\not\sqleq y$, then $g(x)\not\leq g(y)$, so $g(x)_i>g(y)_i$ for some $i$, and $y$ precedes $x$ in $L_i$. Hence $\sqleq$ is the intersection of $L_1$,~\dots,~$L_k$.
\end{proof}

We introduce the following notion, and show how it is related to the concept of dimension:
\begin{definition}
  Given a poset $(P,\sqleq)$, a \emph{sweep system of order $k$} is a family (that is, a $k$-tuple, in which repetitions are allowed) of $k$ linear orders $(\preceq_1,\dots,\preceq_k)$ on $P$ such that, for every interval $I \in \I(P,\sqleq)$ and every $y \in P \setminus I$, there exists an $i$ such that $y \preceq_i x$ for every $x \in I$. The \emph{sweep dimension} of $(P,\sqleq)$, denoted by $\sweepdim(P,\sqleq)$, is the minimum $k$ such that there exists a sweep system of order $k$.
\end{definition}

A sweep system can be seen as a way to identify intervals in a poset by ``sweeping'' the elements of the poset with linear orders. The sweep dimension is a measure of how many linear orders are needed to ensure that every interval can be distinguished from elements outside of it (compare~\cite[Lemma~5.8]{EdJTCG}).
In the following theorem we write $[t\..\infty]$ for the half-line $\{\?x\in\R\mid x\geq t\?\}$ (so that a product of $k$ such half-lines is an upper orthant of $\R^k$), extending the notation for intervals of $[n]$. More precisely:

\begin{theorem}
  \label{thm:orthant}
  Let $(P,\sqleq)$ be a finite poset with $|P|=n$. For every positive integer $k$, the following statements are equivalent:
  \begin{enumerate}
    \item\label{enu:sweepdim} $\sweepdim(P,\sqleq) \leq k$;
    \item\label{enu:sweepmap} there exists a function $\varphi: P \to \R^k$ such that for every interval $I \in \I(P,\sqleq)$ there are $t_1,\dots,t_k$ such that
    \[
      I=\varphi^{-1}\left([t_1\..\infty] \times \dots \times [t_k\..\infty]\right).
    \]
  \end{enumerate}
  Furthermore, the function $\varphi$ of the second item can be chosen so that:
  \begin{itemize}
    \item $\varphi$ takes values in $[n]^k$, each component $\varphi_i: P \to [n]$ is injective, and all the $t_i$'s belong to $[n+1]$ (in particular, the empty interval is the counterimage of the orthant with $t_1=\dots=t_k=n$);
    \item extending $\varphi$ to nonempty intervals by letting
    \[
      \varphi(I)=\Bigl(\min\{\varphi_1(x) \mid x \in I\},\dots,\min\{\varphi_k(x) \mid x \in I\}\Bigr),
    \]
    the thresholds can be chosen so that 
    \[
      I=\varphi^{-1}\left([\varphi(I)_1\..\infty] \times \dots \times [\varphi(I)_k\..\infty]\right);
    \]
    \item for all nonempty intervals $I,J \in \I(P,\sqleq)$, we have $I\supseteq J$ if and only if $\varphi(I) \leq \varphi(J)$ (componentwise); in other words, $I \mapsto \varphi(I)$ is an order embedding of the nonempty intervals, under reverse inclusion, into $([n]^k,\leq)$.
  \end{itemize}
\end{theorem}
\begin{proof}
  \ref{enu:sweepdim} $\implies$ \ref{enu:sweepmap}. Let $|P|=n$, and suppose for the moment $\sweepdim(P,\sqleq) = k$; let $(\preceq_1,\dots,\preceq_k)$ be a sweep system of order $k$. For each $i \in \{1,\dots,k\}$, consider $n$ distinct points in $\R$ and choose a function $\varphi_i: P \to \R$ that maps the elements of $P$ to these points in such a way that the linear order $\preceq_i$ is respected (i.e., if $x \preceq_i y$ then $\varphi_i(x) \leq \varphi_i(y)$). Define the function $\varphi: P \to \R^k$ by letting $\varphi(x)=(\varphi_1(x),\dots,\varphi_k(x))$. 

  Let $I$ be an interval of $(P,\sqleq)$. If $I$ is empty, we can take each $t_i$ to be larger than the maximum value of $\varphi_i$, so that $I=\varphi^{-1}\left([t_1\..\infty] \times \dots \times [t_k\..\infty]\right)$. If $I$ is nonempty, for every $i\in \{1,\dots,k\}$ let $t_i=\min \{\varphi_i(x) \mid x \in I\}$. Then, for every $x \in I$, we have $\varphi(x)_i =\varphi_i(x) \geq t_i$, so
  \[
    I \subseteq \varphi^{-1}\left([t_1\..\infty] \times \dots \times [t_k\..\infty]\right).
  \]
  To prove the other inclusion, let $y \in P \setminus I$. By the definition of sweep system, there exists an $i \in \{1,\dots,k\}$ such that $y \preceq_i x$ for every $x \in I$, that is, $\varphi(y)_i < t_i$. So $y \notin \varphi^{-1}\left([t_1\..\infty] \times \dots \times [t_k\..\infty]\right)$.

  Repetitions are allowed in a sweep system, so the case $\sweepdim(P,\sqleq)<k$ is covered by repeating sufficiently many times one of the linear orders of a sweep system of minimum order.

  \ref{enu:sweepmap} $\implies$ \ref{enu:sweepdim}. Define $\preceq_i$ as the linear order on $P$ induced by the $i$-th component of $\varphi$, that is, $x \preceq_i y$ iff $\varphi(x)_i \leq \varphi(y)_i$, with ties broken arbitrarily. Let $I \in \I(P,\sqleq)$ and let $y \in P \setminus I$. By hypothesis, there exist $t_1,\dots,t_k$ such that
  \[
    I=\varphi^{-1}\left([t_1\..\infty] \times \dots \times [t_k\..\infty]\right).
  \]
  Since $y \notin I$, there exists an $i$ such that $\varphi(y)_i < t_i$, that is, $y \preceq_i x$ for every $x \in I$. Hence $(\preceq_1,\dots,\preceq_k)$ is a sweep system of order $k$, so $\sweepdim(P,\sqleq) \leq k$.

  \smallskip
  Finally, let us prove the last three items. The first one is immediate from the construction of $\varphi$ in the first part of the proof: each $\varphi_i$ maps $P$ bijectively onto $n$ distinct points in $\R$, which we may take to be the elements of $[n]$; moreover the trace of the orthant with $t_1=\dots=t_k=n$ is empty, as $\varphi$ takes values in $[n]^k$.

  For the second item, let $I$ be a nonempty interval and let $t_1,\dots,t_k$ be such that $I=\varphi^{-1}\left([t_1\..\infty] \times \dots \times [t_k\..\infty]\right)$. Then $t_i \leq \varphi_i(x)$ for every $x \in I$, so $t_i \leq \varphi(I)_i$; on the other hand $\varphi_i(x)\geq \varphi(I)_i$ for every $x \in I$ by definition of $\varphi(I)$. Hence
  \begin{multline*}
    I \subseteq \varphi^{-1}\left([\varphi(I)_1\..\infty] \times \dots \times [\varphi(I)_k\..\infty]\right) \subseteq \\ \subseteq \varphi^{-1}\left([t_1\..\infty] \times \dots \times [t_k\..\infty]\right)=I,
  \end{multline*}
  and all inclusions are equalities.

  For the third item: if $I \supseteq J$, then for every $i$ we have $\varphi(I)_i = \min\{\varphi_i(x) \mid x \in I\} \leq \min\{\varphi_i(x) \mid x \in J\} =  \varphi(J)_i$, because a minimum computed over a superset can only be smaller: so $\varphi(I) \leq \varphi(J)$. Conversely, if $\varphi(I) \leq \varphi(J)$, then, using the second item twice,
  \begin{multline*}
    J \subseteq \varphi^{-1}\left([\varphi(J)_1\..\infty] \times \dots \times [\varphi(J)_k\..\infty]\right) \subseteq\\ \subseteq \varphi^{-1}\left([\varphi(I)_1\..\infty] \times \dots \times [\varphi(I)_k\..\infty]\right)=I,
  \end{multline*}
  where the first containment holds because $\varphi_i(x) \geq \varphi(J)_i$ for every $x \in J$, and the second because $\varphi(I)\leq \varphi(J)$ makes the first orthant larger than the second.
\end{proof}

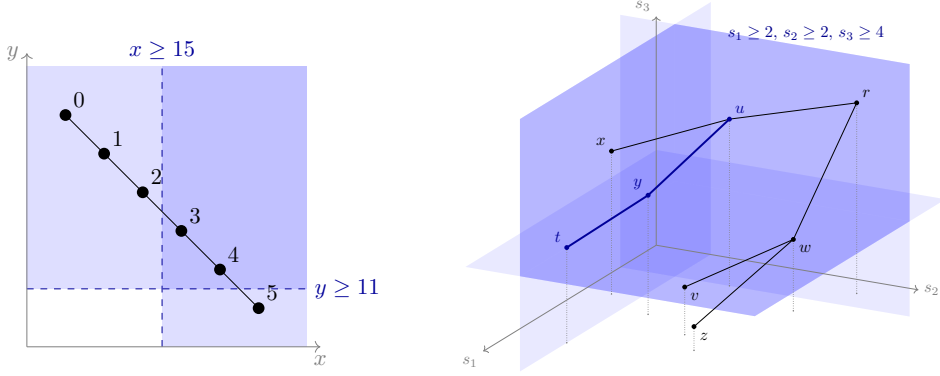
\begin{figure}
  \centering
  \resizebox{0.4\textwidth}{!}{\begin{tikzpicture}[x=.3cm,y=.3cm]
  \fill[blue,opacity=.13] (15,8) rectangle (22.5,22.5);
  \fill[blue,opacity=.13] (8,11) rectangle (22.5,22.5);
  \draw[blue!60!black,dashed] (15,8) -- (15,22.5) node[above]{$x\geq 15$};
  \draw[blue!60!black,dashed] (8,11) -- (22.5,11) node[right]{$y\geq 11$};
  \draw[->,gray] (8,8) -- (23.2,8) node[below]{$x$};
  \draw[->,gray] (8,8) -- (8,23.2) node[left]{$y$};
  \draw[thin] (10,20) -- (12,18) -- (14,16) -- (16,14) -- (18,12) -- (20,10);
  \foreach \x/\y/\i in {10/20/0,12/18/1,14/16/2,16/14/3,18/12/4,20/10/5}{
    \fill (\x,\y) circle (2.6pt);
    \node[above right=0pt] at (\x,\y) {$\i$};
  }
\end{tikzpicture}}$\qquad$ 
  \resizebox{0.5\textwidth}{!}{\begin{tikzpicture}[x={(-0.62cm,-0.38cm)},y={(0.94cm,-0.16cm)},z={(0cm,0.82cm)},scale=.62]
  \fill[blue,opacity=.09] (2,0,0) -- (2,10.6,0) -- (2,10.6,10.6) -- (2,0,10.6) -- cycle;
  \fill[blue,opacity=.09] (0,2,0) -- (10.6,2,0) -- (10.6,2,10.6) -- (0,2,10.6) -- cycle;
  \fill[blue,opacity=.09] (0,0,4) -- (10.6,0,4) -- (10.6,10.6,4) -- (0,10.6,4) -- cycle;
  \fill[blue,opacity=.14] (2,2,4) -- (2,10.6,4) -- (2,10.6,10.6) -- (2,2,10.6) -- cycle;
  \fill[blue,opacity=.14] (2,2,4) -- (10.6,2,4) -- (10.6,2,10.6) -- (2,2,10.6) -- cycle;
  \fill[blue,opacity=.14] (2,2,4) -- (10.6,2,4) -- (10.6,10.6,4) -- (2,10.6,4) -- cycle;
  \fill[blue,opacity=.09] (10.6,2,4) -- (10.6,10.6,4) -- (2,10.6,4)
    -- (2,10.6,10.6) -- (2,2,10.6) -- (10.6,2,10.6) -- cycle;
  \draw[->,gray] (0,0,0) -- (9.6,0,0) node[below left]{$s_1$};
  \draw[->,gray] (0,0,0) -- (0,9.6,0) node[right]{$s_2$};
  \draw[->,gray] (0,0,0) -- (0,0,9.6) node[above left]{$s_3$};
  \node[blue!55!black,above] at (2,6.8,10.6) {$s_1\geq2$,\ $s_2\geq2$,\ $s_3\geq4$};
  \foreach \px/\py/\pz in {1/8/8, 2/4/7, 3/7/3, 4/1/6, 5/3/5, 6/5/2, 7/6/1, 8/2/4}{
    \draw[densely dotted,gray] (\px,\py,0) -- (\px,\py,\pz);
    \fill[gray] (\px,\py,0) circle (.7pt);}
  \draw[thin] (1,8,8) -- (2,4,7)  (1,8,8) -- (3,7,3)  (2,4,7) -- (4,1,6)
              (3,7,3) -- (6,5,2)  (3,7,3) -- (7,6,1);
  \draw[very thick,blue!60!black] (2,4,7) -- (5,3,5) -- (8,2,4);
  \fill[blue!60!black] (2,4,7) circle (2.4pt) node[above right=0pt]{$u$};
  \fill[blue!60!black] (5,3,5) circle (2.4pt) node[above left=0pt]{$y$};
  \fill[blue!60!black] (8,2,4) circle (2.4pt) node[above left=0pt]{$t$};
  \fill (1,8,8) circle (2.4pt) node[above right=0pt]{$r$};
  \fill (3,7,3) circle (2.4pt) node[below right=0pt]{$w$};
  \fill (4,1,6) circle (2.4pt) node[above left=0pt]{$x$};
  \fill (6,5,2) circle (2.4pt) node[below right=0pt]{$v$};
  \fill (7,6,1) circle (2.4pt) node[below right=0pt]{$z$};
\end{tikzpicture}}
  \caption{\label{fig:orthant}Left: a sweep system of order 2 for a linear order (singling out the interval $[3\..4]$). Right: three sweeps identifying the intervals of the form $[a\..b]$ (the vertical paths) of a tree, here singling out $[u\..t]$; identifying \emph{all} intervals of this tree would require a sweep system of order at least five (see text). The tree has root $r$, and its cover relations are $r\sqless u$, $r\sqless w$, $u\sqless x$, $u\sqless y$, $y\sqless t$, $w\sqless v$ and $w\sqless z$. In both cases, the partial order is identified through its Hasse diagram.}
\end{figure}

The previous result gives a geometrical interpretation to a sweep system of order $k$ for $(P,\sqleq)$: the existence of such a system makes it possible to represent the elements of $P$ as points in the Euclidean space $\R^k$ so that each interval can be identified by an upper orthant. Figure~\ref{fig:orthant} (left) shows an example: a sweep system of order 2 for the linear order $([6],\leq)$. In the picture, we identify the interval $[3\..4]$ by two half-planes ($x\geq 15$ and $y \geq 11$); moreover, in this case any two similar half-planes identify a (possibly empty) interval.

In Figure~\ref{fig:orthant} (right) we consider a poset $(P,\sqleq)$ that is a tree: a \emph{tree} is a poset such that, for every $x \in P$, the set $\{y \in P \mid y \sqleq x\}$ is linearly ordered by $\sqleq$ (strictly speaking this defines a forest; the tree of the figure has a minimum). Here, however, three sweeps cannot identify \emph{all} intervals: every subset of $\{\?w,x,y,v,z\?\}$ is an interval (as one can easily check on the figure: the only comparable pairs among these five elements are $w\sqless v$ and $w\sqless z$, and no element of the tree lies strictly between $w$ and $v$ or between $w$ and $z$), and this forces every sweep system of this tree to have order at least five, as we shall see (Proposition~\ref{prop:sweepprops}). What the figure shows is that three sweeps suffice to identify all the intervals \emph{of the form} $[a\..b]=\{\?c \in P \mid a \sqleq c \sqleq b\?\}$: the tree is embedded in $\R^3$, and the orthant shown identifies $[u\..t]$.

\subsection{Staircase regions and sweep dimension}

An \emph{$n$-staircase region of dimension $k$} is a subset $S$ of $[n]^k$ such that for every ${\mathbf x} \in S$ and every ${\mathbf y} \in [n]^k$, if ${\mathbf y} \leq {\mathbf x}$ (componentwise), we have ${\mathbf y} \in S$: that is, $S$ is a lower set of $([n]^k,\leq)$ (equivalently, its complement is a union of upper orthants). The set of all $n$-staircase regions of dimension $k$ is denoted by $\Stair_{n,k}$.

Up to a coordinate change, $n$-staircase functions correspond to the $n$-staircase regions of dimension $2$ contained in the triangle $T_n$ of Section~\ref{sec:partitions}, the correspondence reversing the order (pointwise order on functions, inclusion on regions). More precisely (through Theorem~\ref{thm:delta}, this is the correspondence between Dyck paths and partitions contained in a staircase of~\cite[Proposition~3.9]{FeMLPRTV} and~\cite[Figure~2]{FerUDP}, restated for staircase functions):

\begin{lemma}
  \label{lem:triangle}
  Let $\mathcal T_n$ be the set of $n$-staircase regions of dimension $2$ contained in the triangle $T_n$. For $f\in\Stair_n$, recall from Section~\ref{sec:partitions} that
  \[
    S_f=\{\,(x,b)\in[n]^2 \mid b+f(x)\leq n-1\,\}.
  \]
  Then $f\mapsto S_f$ is a bijection between $\Stair_n$ and $\mathcal T_n$, whose inverse maps $S$ to $f_S(x)=n-|\{\,b \mid (x,b)\in S\,\}|$; moreover $f\leq g$ iff $S_f\supseteq S_g$.
\end{lemma}
\begin{proof}
  $S_f$ is a lower set: if $(x,b)\in S_f$ and $(x',b')\leq(x,b)$, then $b'+f(x')\leq b+f(x)\leq n-1$ because $f$ is non-decreasing. It is contained in $T_n$ because $f(x)\geq x$. The column of $S_f$ over $x$ is $\{\,b \mid b\leq n-1-f(x)\,\}$, which has $n-f(x)$ elements, so $f_{S_f}=f$ and the map is injective.

  Conversely, let $S\in\mathcal T_n$ and let $h_x=|\{\,b \mid (x,b)\in S\,\}|$ be the height of the column over $x$. Since $S$ is a lower set, each column is an initial segment $\{\,0,\dots,h_x-1\,\}$ of $[n]$, and the column over $x'\leq x$ contains the column over $x$, so $x\mapsto h_x$ is non-increasing; since $S\subseteq T_n$, every $(x,b)\in S$ has $b\leq n-1-x$, so $h_x\leq n-x$. Hence $f_S(x)=n-h_x$ satisfies $x\leq f_S(x)\leq n$ and is non-decreasing, that is, $f_S\in\Stair_n$, and $S_{f_S}=\{\,(x,b) \mid b\leq h_x-1\,\}=S$. So the map is surjective.

  Finally $f\leq g$ iff $n-f(x)\geq n-g(x)$ for all $x$, that is, iff every column of $S_f$ contains the corresponding column of $S_g$, that is, iff $S_f\supseteq S_g$.
\end{proof}

We can now generalize to an arbitrary finite poset the correspondence between antichains and staircase functions (that is, by Theorem~\ref{thm:delta}, Dyck words) provided by Theorem~\ref{thm:rho}. Fix a sweep system of order $k$ for $(P,\sqleq)$, and let $\varphi$ be a function as in the second part of Theorem~\ref{thm:orthant}.

\begin{definition}
  For $A \in \E(P,\sqleq)\setminus\{\top\}$, the \emph{staircase region} of $A$ is the lower set generated in $[n]^k$ by the images of its witnesses:
  \[
    R(A) = \{\,{\mathbf t} \in [n]^k \mid {\mathbf t} \leq \varphi(I) \text{ for some } I \in A\,\}.
  \]
\end{definition}

Note that the definition makes sense: since $A \neq \top$, all the witnesses of $A$ are nonempty (an antichain containing $\emptyset$ can only be $\{\emptyset\}=\top$, as the empty interval is contained in every interval), so $\varphi(I)$ is defined for every $I \in A$.

\begin{theorem}
  \label{thm:region}
  The map $A \mapsto R(A)$ is an order embedding of $\E(P,\sqleq)\setminus\{\top\}$ into $(\Stair_{n,k},\subseteq)$, that is, $A \leq B$ iff $R(A) \subseteq R(B)$. Moreover:
  \begin{enumerate}
    \item\label{enu:reg-join} $R(A \OR B)=R(A) \cup R(B)$;
    \item\label{enu:reg-max} the maximal points of $R(A)$ are exactly the images $\varphi(I)$ of the witnesses $I \in A$; in particular, $A$ can be recovered from $R(A)$.
  \end{enumerate}
\end{theorem}
\begin{proof}
  Suppose $A \leq B$, and let $I \in A$: there is $J \in B$ with $J \subseteq I$, and by the last item of Theorem~\ref{thm:orthant} $\varphi(I) \leq \varphi(J)$, so every generator of $R(A)$ belongs to $R(B)$, hence $R(A) \subseteq R(B)$. Conversely, suppose $R(A) \subseteq R(B)$, and let $I \in A$: since $\varphi(I) \in R(B)$, there is $J \in B$ with $\varphi(I) \leq \varphi(J)$, which by the last item of Theorem~\ref{thm:orthant} means $J \subseteq I$: hence $A \leq B$.

  For item~\ref{enu:reg-join}, recall~\cite{BoVLAFI} that $A \OR B$ is the set of $\subseteq$-minimal elements of $A \cup B$; if a witness $I \in A \cup B$ is discarded because $J \subseteq I$ for some $J \in A \OR B$, then $\varphi(I) \leq \varphi(J)$, so discarded witnesses generate no new points: the lower set generated by the images of $A \OR B$ is the same as the one generated by the images of $A \cup B$, which is $R(A) \cup R(B)$.

  For item~\ref{enu:reg-max}, the images of distinct witnesses of $A$ are incomparable (by the last item of Theorem~\ref{thm:orthant}, comparable images would mean one witness contained in another, impossible in an antichain), so each $\varphi(I)$ is a maximal point of $R(A)$, and no other point of $R(A)$ is maximal (every point lies below some $\varphi(I)$ by definition). Since $I \mapsto \varphi(I)$ is injective (being an order embedding), $A$ is determined by $R(A)$.
\end{proof}

\begin{remark}
  \label{rem:lowersets}
  Theorem~\ref{thm:region} is an instance of a general principle. Let $Q=(\I(P,\sqleq)\setminus\{\emptyset\},\supseteq)$: $\E(P,\sqleq)\setminus\{\top\}$ is the lattice of lower sets of $Q$ (a lower set being generated by the antichain of its maximal elements, that is, by the $\subseteq$-minimal intervals it contains), with join given by union. An order embedding of $Q$ into any poset $M$ induces an embedding of the lattice of lower sets of $Q$ into that of $M$, mapping a lower set to the lower set it generates, and preserving unions; the last item of Theorem~\ref{thm:orthant} says that $\varphi$ is such an embedding of $Q$ into $M=([n]^k,\leq)$. What the sweep formulation adds is that $\varphi$ is read off $k$ linear orders of $P$, that the maximal points of $R(A)$ are the images of the witnesses, and that for $k=2$ and $P$ a chain it recovers Dyck words (Proposition~\ref{prop:regionchain}).
\end{remark}

The correspondence of Theorem~\ref{thm:rho} is precisely the case $k=2$:

\begin{proposition}
  \label{prop:regionchain}
  Let $P=[n]$ with the usual order, and consider the sweep system $\{\leq,\geq\}$ with $\varphi_1(x)=x$ and $\varphi_2(x)=n-1-x$. Then $\varphi([x\..y])=(x,n-1-y)$ for every nonempty interval, and for every $A \in \E_n\setminus\{\top\}$
  \[
    R(A) = S_{\rho_A}.
  \]
  Moreover, in this case the map $A \mapsto R(A)$ is a bijection of $\E_n\setminus\{\top\}$ onto the set $\mathcal T_n$ of staircase regions contained in the triangle.
\end{proposition}
\begin{proof}
  The first statement is trivial. For the second, the column of $R(A)$ over $a \in [n]$ consists of the points $(a,b)$ with $b \leq n-1-y$ for some $[x\..y] \in A$ with $x \geq a$: the largest such $b$ is $n-1-\rho_A(a)$ by definition of $\rho_A$ (if no witness $[x\..y]\in A$ has $x\geq a$, then $\rho_A(a)=n$ and both columns are empty), so the column coincides with the column of $S_{\rho_A}$. For the last statement, $A\mapsto R(A)=S_{\rho_A}$ is the composition of the bijection $A\mapsto\rho_A$ of $\E_n\setminus\{\top\}$ onto $\Stair_n$ (Theorem~\ref{thm:rho}) with the bijection $f\mapsto S_f$ of $\Stair_n$ onto $\mathcal T_n$ (Lemma~\ref{lem:triangle}).
\end{proof}

Summing up: for a chain, the elements of the antichain lattice \emph{are} the two-dimensional staircase regions of a triangle (equivalently: staircase functions,  Dyck words) with joins given by unions; for a general poset $(P,\sqleq)$, they embed, join-compatibly, into the $k$-dimensional staircase regions, and by Theorem~\ref{thm:orthant} the sweep dimension is exactly the least $k$ for which the elements of $P$ can be placed in $[n]^k$ so that every interval is cut out by an upper orthant, which is what the construction of $R(A)$ requires. The objects replacing Dyck words in higher sweep dimension are thus the natural higher-dimensional staircases: for $k=3$, stacks of unit cubes pushed into a corner, and in general regions bounded by a monotone $(k-1)$-dimensional staircase hypersurface; equivalently, slicing along any one coordinate, a region of dimension $k$ is a weakly nested family of $n$ regions of dimension $k-1$ (for $k=3$: $n$ nested staircase paths). Two caveats distinguish the general case from the chain case, and both are unavoidable: the map $A\mapsto R(A)$ is an embedding but not, in general, a bijection onto $\Stair_{n,k}$ (which regions arise depends on $(P,\sqleq)$ and on the chosen sweep system: this must be so, since, for instance, the number of elements of $\Stair_{n,3}$ is a plane-partition count, while $|\E(P,\sqleq)|$ varies with $P$); and only joins are computed by unions, while the meet of $\E(P,\sqleq)$ is in general finer than the intersection of the corresponding regions.

\subsection{Properties of sweep dimension}

We have the following:
\begin{proposition}
  \label{prop:dimsweep}
  For any poset $(P,\sqleq)$ with $|P|\geq1$ we have that
  \[
    \dim\Bigl(\I(P,\sqleq), \supseteq\Bigr) = \sweepdim(P,\sqleq).
  \]
\end{proposition}
\begin{proof}
  We start by showing the $\geq$ inequality. Recall (Lemma~\ref{lem:dimembed}) that a poset has dimension at most $k$ if and only if it can be order-embedded in $\R^k$ with coordinatewise order. Let $k=\dim\Bigl(\I(P,\sqleq),\supseteq\Bigr)$ and $f: \I(P,\sqleq) \to \R^k$ be an order embedding. For every $x \in P$, let ${\mathbf a}_x=f([x\..x]) \in \R^k$. Note that
  $I \supseteq [x\..x]$ iff $x \in I$, so for every interval $I$ and every $x \in I$, we must have $f(I)\leq {\mathbf a}_x$ (componentwise), so for every $i$,
  \[
    \min \{({\mathbf a}_x)_i \mid x \in I\}\geq (f(I))_i.
  \]
  On the other hand, if $y \not\in I$, then $I \not\supseteq [y\..y]$, hence necessarily $f(I)\not\leq {\mathbf a}_y$. Hence there must be an index $i \in \{1,\dots,k\}$ such that $f(I)_i > ({\mathbf a}_y)_i$, so
  \[
    \min \{({\mathbf a}_x)_i \mid x \in I\} \geq (f(I))_i > ({\mathbf a}_y)_i.
  \]
  So the $k$ coordinatewise orders $\preceq_1$, \dots, $\preceq_k$ defined on $P$ by 
  \[ 
  y \preceq_i x \text{ iff } ({\mathbf a}_y)_i \leq ({\mathbf a}_x)_i
  \]
  (ties broken arbitrarily) are a sweep system of order $k$ for $(P, \sqleq)$, hence $k \geq \sweepdim(P,\sqleq)$.

  For the other inequality, let $\sweepdim(P,\sqleq)=k$ and $(\preceq_1,\dots,\preceq_k)$ be a sweep system; also let $|P|=n$. For every $i$, let $p_i(x) \in \{1,\dots,n\}$ be the rank of $x \in P$ in $\preceq_i$, and for every interval $I$ let $m_i(I)=\min \{p_i(x) \mid x \in I\}$, with the convention $m_i(\emptyset)=n+1$ (so that $m_i(I)\leq n<m_i(\emptyset)$ for every nonempty $I$). Define $g_i: \I(P,\sqleq) \to \R$ by letting
  \[
    g_i(I)=(n+1)\, m_i(I) - |I|.
  \]
  If $I \supseteq J$, then
  \[
    m_i(I) \leq m_i(J)
  \]
  (this holds also when $J=\emptyset$, by the convention above), so either the two minima are equal (hence $g_i(I) \leq g_i(J)$, since $|I|\geq |J|$), or
  \[
    g_i(J)-g_i(I) \geq n+1+|I|-|J|\geq n+1>0,
  \]
  which again implies $g_i(I) \leq g_i(J)$.

  On the other hand, if $I \not\supseteq J$, let $y \in J \setminus I$ and consider an index $i$ such that $y \preceq_i x$ for all $x \in I$. Then
  \[
    m_i(J) \leq p_i(y) \leq m_i(I)-1
  \]
  (the second inequality holds also when $I=\emptyset$, as $p_i(y)\leq n$), so
  \begin{multline*}
    g_i(J) \leq (n+1) \bigl(m_i(I)-1\bigr) - |J|=\\
    =(n+1)\, m_i(I) - (n+1) - |J| <\\
    <(n+1)\, m_i(I) - |I| = g_i(I),
  \end{multline*}
  because $|I|\leq n<n+1+|J|$; so $g_i(I) \not\leq g_i(J)$.

  We have proven that
  \[
      I \supseteq J \qquad\text{ iff }\qquad g_i(I) \leq g_i(J) \quad \forall i \in \{1,\dots,k\}.
  \]
  Hence the function $g: \I(P,\sqleq) \to \R^k$ defined by $g(I)=(g_1(I),\dots,g_k(I))$ is an order embedding, hence the dimension of $\I(P,\sqleq)$ is at most $k$.
\end{proof}

In other words, the sweep dimension of a poset $P$ measures the standard dimension of its interval poset $\I(P)$. Two consequences are worth recording. First, $\dim(P,\sqleq)\leq\sweepdim(P,\sqleq)$: the map $x\mapsto\downset x$ is an order embedding of $(P,\sqleq)$ into $(\I(P,\sqleq),\subseteq)$, which has the same dimension as its dual $(\I(P,\sqleq),\supseteq)$, and the dimension of a poset is at least that of its subposets. Second, the lattice $\E(P,\sqleq)\setminus\{\top\}$ is the lattice of lower sets of the poset $Q$ of Remark~\ref{rem:lowersets}, and the dimension of the lattice $L$ of lower sets of a finite nonempty poset $Q$ is $w=\width(Q)$. Indeed, if $C_1$,~\dots,~$C_w$ is a chain decomposition of $Q$, the map $D\mapsto(|D\cap C_1|,\dots,|D\cap C_w|)$ is an order embedding of $L$ into $\R^w$, as each $D\cap C_i$ is an initial segment of $C_i$, so $\dim L\leq w$ by Lemma~\ref{lem:dimembed} (compare Dilworth's embedding theorem for distributive lattices~\cite[Theorem~1.2]{DilDTPOS}); conversely, if $A$ is an antichain of $Q$ with $w$ elements, the lower sets $(\downset A\setminus A)\cup\{a\}$ and $\downset A\setminus\{a\}$, for $a\in A$, form a copy of the standard example of Dushnik and Miller, which has dimension $w$~\cite[Theorem~4.1]{DuMPOS}. Hence $\dim(\E(P,\sqleq)\setminus\{\top\})$ is the size of the largest family of pairwise $\subseteq$-incomparable nonempty intervals of $P$; for $P=[n]$ this gives $n$, the dimension of the Stanley lattice of semilength $n+1$. 
The following proposition collects the basic properties of the sweep dimension, and in particular it connects it to the width of a poset.

\begin{proposition}
  \label{prop:sweepprops}
  Let $(P,\sqleq)$ be a poset with $|P|\geq 2$, and let $\theta(P,\sqleq)$ be the maximum cardinality of a subposet of $(P,\sqleq)$ with no $3$-element chain. Then:
  \begin{enumerate}
    \item\label{enu:sd-dual} $\sweepdim(P,\sqleq) = \sweepdim(P,\sqgeq)$;
    \item\label{enu:sd-mono} if $P' \subseteq P$, then $\sweepdim(P',\sqleq') \leq \sweepdim(P,\sqleq)$ (where $\sqleq'$ is the restriction on $P'$ of $\sqleq$);
    \item\label{enu:sd-twolayer} $\theta(P,\sqleq) \leq \sweepdim(P,\sqleq)$;
    \item\label{enu:sd-bounds} $2 \leq \sweepdim(P,\sqleq) \leq 2\cdot\width(P,\sqleq)$.
  \end{enumerate}
\end{proposition}
Items~\ref{enu:sd-twolayer} and~\ref{enu:sd-bounds} can also be derived from~\cite[Theorems~4.1 and~4.6, Proposition~4.4]{EdSCRCDCG}, as explained in Section~\ref{sec:related}.
\begin{proof}
  Item~\ref{enu:sd-dual} is a trivial consequence of Proposition~\ref{prop:dimsweep} (intervals of a poset and of its opposite are the same).

  For item~\ref{enu:sd-mono}, note first that, for every $S \subseteq P$, the set
  \[
    \overline S = \{\,z \in P \mid s \sqleq z \sqleq t \text{ for some } s,t \in S\,\}
  \]
  is an interval of $(P,\sqleq)$ containing $S$: if $z_1,z_2 \in \overline S$ (say $s_1 \sqleq z_1 \sqleq t_1$ and $s_2 \sqleq z_2 \sqleq t_2$ with $s_i,t_i \in S$), and $z_1 \sqleq w \sqleq z_2$, then $s_1 \sqleq w \sqleq t_2$, so $w \in \overline S$. Now suppose that $(\preceq_1,\dots,\preceq_k)$ is a sweep system for $(P,\sqleq)$, and let $I' \in \I(P',\sqleq')$ and $y \in P' \setminus I'$. Every $z \in \overline{I'} \cap P'$ lies between two elements of $I'$, hence $z \in I'$ by convexity of $I'$ in $(P',\sqleq')$: as a consequence, in particular, $y \notin \overline{I'}$. Let $i$ be such that $y \preceq_i x$ for every $x \in \overline{I'}$: since $I' \subseteq \overline{I'}$, the restriction of $\preceq_i$ to $P'$ lists $y$ before every element of $I'$. Hence the restrictions of $\preceq_1,\dots,\preceq_k$ to $P'$ form a sweep system for $(P',\sqleq')$.

  For item~\ref{enu:sd-twolayer}, by item~\ref{enu:sd-mono} it suffices to show that $\sweepdim(P',\sqleq')\geq|P'|$ when the induced subposet contains no $3$-element chain. In this case no element of $P'$ lies strictly between two others, so the convexity condition is vacuous and \emph{every} nonempty subset of $P'$ is an interval. Given a sweep system, for every $y \in P'$, $y$ must be smaller than all the elements of the set (i.e., interval) $P'\setminus\{y\}$ in some sweep; since distinct elements require distinct sweeps, there must be at least $|P'|$ of them.

  For item~\ref{enu:sd-bounds}, the lower bound is the case $|P'|=2$ of item~\ref{enu:sd-twolayer} (any two distinct elements induce a subposet without $3$-element chains). For the upper bound, let $w = \width(P,\sqleq)$: using Dilworth's theorem, cover $P$ with chains $C_1,\dots,C_w$. For each $j$, consider the two sweeps listing first the elements of $C_j$ (in ascending, respectively descending, order with respect to $\sqleq$) and then the remaining elements of $P$, in any order. Let now $I$ be an interval and $y \not\in I$; assume $y \in C_j$. The set $I \cap C_j$ is an interval of the subposet induced on the chain $C_j$; since $y \in C_j$ but $y \not\in I \cap C_j$, $y$ lies strictly below all of it or strictly above all of it in $C_j$. Hence, one of the two sweeps corresponding to $C_j$ lists $y$ before every element of $C_j \cap I$ and (listing $C_j$ first) before every other element of $I$ as well.
\end{proof}

Note that the sweep dimension is governed by the width of the poset, not by its dimension: an antichain with $n$ elements has sweep dimension $n$, although its dimension is $2$. In particular, the sweep dimension is not bounded by any function of the dimension. Item~\ref{enu:sd-twolayer} also justifies the claim made about the tree of Figure~\ref{fig:orthant}: the five elements $w,x,y,v,z$ induce a subposet with no $3$-element chain, so the sweep dimension of that tree is at least five.

A consequence of the previous proposition is that:
\begin{corollary}
  \label{cor:sweepprops}
  Let $(P,\sqleq)$ be a poset with $|P|\geq 2$. Then:
  \[
    \max(2,\width(P,\sqleq)) \leq \sweepdim(P,\sqleq) \leq 2\cdot\width(P,\sqleq).
  \]
\end{corollary}
\begin{proof}
  The lower bound is a consequence of Proposition~\ref{prop:sweepprops} (\ref{enu:sd-twolayer}) and (\ref{enu:sd-bounds}), observing that $\theta(P,\sqleq) \geq \width(P,\sqleq)$.
\end{proof}

The two bounds of Corollary~\ref{cor:sweepprops} are tight: 

\begin{proposition}
  For every $w$ and $s$ with $\max(2,w) \leq s \leq 2w$, there exists a poset $(P,\sqleq)$ such that $\width(P,\sqleq)=w$ and $\sweepdim(P,\sqleq)=s$.
\end{proposition}
\begin{proof}
  If $w=1$, then $s=2$ and the poset $(P,\sqleq)$ with $P=\{x,y\}$ and $x \sqleq y$ has width $1$ and sweep dimension $2$. If $1<w$, an antichain of size $w$ has sweep dimension $s=w$ (it must be at least $w$ by Proposition~\ref{prop:sweepprops} (\ref{enu:sd-twolayer}), and clearly a sweep system of order $w$ can be obtained by listing each element first in one of the sweeps).

  We are left to cover the case $1<w<s\leq 2w$.
  Consider a complete bipartite poset $(P,\sqleq)$ with $w$ minimal elements and $s-w$ maximal elements, every minimal element being smaller than every maximal element. Then $\width(P,\sqleq)=w$: an antichain cannot contain both a minimal and a maximal element, so it is contained either in the set of minimal elements or in the set of maximal elements, and the former is the larger of the two, because $s-w\leq 2w-w=w$. Moreover $\sweepdim(P,\sqleq)=s$: on the one hand $\sweepdim(P,\sqleq)\leq |P|=s$, because a sweep system of order $|P|$ is obtained by listing each element first in one sweep; on the other hand $(P,\sqleq)$ contains no $3$-element chain, so $\theta(P,\sqleq)=|P|=s$ and $\sweepdim(P,\sqleq)\geq s$ by Proposition~\ref{prop:sweepprops} (\ref{enu:sd-twolayer}).
\end{proof}

\section{Comparison with the literature}
\label{sec:related}

\paragraph{Lattices of Dyck paths.} Dyck paths of semilength $m$ carry three
lattice structures, the \emph{Stanley}, \emph{Tamari}, and \emph{Kreweras}
lattices; the first one is the order in which a Dyck path is below another if it
lies weakly below it, and it is a distributive lattice. The three lattices are
ordered by inclusion: the Stanley lattice is an extension of the Tamari
lattice, which is an extension of the Kreweras lattice~\cite{BeBICLRT}. The
Tamari and Kreweras lattices are not distributive~\cite[p.~222]{StaFL}, so the
Stanley lattice is the only one of the three that can be isomorphic (or
anti-isomorphic) to $\E_n$.

Ferrari and Pinzani~\cite{FePLLP} gave a criterion for posets of lattice
paths ordered by dominance to be distributive lattices, which applies in
particular to Dyck paths (and also to Motzkin and Schr\"oder paths), and
proved that the rank of a Dyck path of semilength $m$ in the Stanley lattice
is half the difference between the sum of its heights and
$m$~\cite[Proposition~3.1]{FePLLP}, which is our area statistic. They also
state that the join-irreducible elements of the Stanley
lattice form the poset of intervals of a chain. The statement was later proved by Ferrari and Munarini~\cite[Proposition~3.7]{FeMLPRTV}: the
join-irreducible elements of the Stanley lattice of semilength $m$ are the
paths with exactly one ``pyramid'', and they form a poset isomorphic to the
poset of intervals of an $(m-1)$-element chain; hence, by Birkhoff's
representation theorem, the Stanley lattice of semilength $m$ is isomorphic to
the lattice of order ideals of this poset. Ferrari and Munarini also show that
the Stanley lattice is anti-isomorphic to the lattice of partitions
fitting inside the staircase $(m-1,m-2,\dots,1)$~\cite[Proposition~3.9]{FeMLPRTV},
and mention that its rank-unimodality is still open.

The description of Dyck paths by nondecreasing sequences compared
coordinatewise, which underlies our staircase functions, is also known.
Sapounakis, Tasoulas and Tsikouras~\cite{SaTTDPODP} encode a Dyck path by
its \emph{dominating sequence} (the numbers of rises between consecutive
falls), so that the dominance order becomes the comparison of partial sums,
which are our quantities $\alpha_k(w)$; they give the rank of a path in these
terms~\cite[equation (1)]{SaTTDPODP}, a recursion for the rank generating function
that is the classical recursion for the Carlitz--Riordan
polynomials~\cite[Proposition~3.2]{SaTTDPODP}, and the M\"obius function of
the lattice; they also credit Narayana and Fulton with the distributivity of
the dominance order on Grand-Dyck paths, of which the Stanley lattice is a
sublattice. Barcucci, Bernini, Ferrari and Poneti~\cite{BBFPDLS} observe
that the increasing parking functions (nondecreasing sequences with
$a_i\leq i$), ordered coordinatewise, form a lattice anti-isomorphic to
the Stanley lattice, and prove that the Stanley lattice is isomorphic to the
lattice of $312$-avoiding permutations under the strong Bruhat
order~\cite[Corollary~5.1]{BBFPDLS}. Finally, Ferrari and
Munarini~\cite{FeMECSCDL} enumerate the chains of every length and the
saturated chains of length at most four of the Stanley lattice, results that
transfer verbatim to $\E_n\setminus\{\top\}$.

Earlier,
Stanley's paper on the Fibonacci lattice~\cite{StaFL} had observed that the
order ideals of the poset of intervals of an $n$-element chain, ordered by
inclusion, are counted by the Catalan number $C_{n+1}$, and described a
bijection with the lattice paths that never go below the diagonal.

Since $\E_n\setminus\{\top\}$ is anti-isomorphic to the lattice of these
order ideals, the lattice-theoretic content of our anti-isomorphism between
$\E_n\setminus\{\top\}$ and the Stanley lattice follows from the
results of Ferrari and Munarini through Birkhoff's representation theorem,
although no explicit isomorphism with antichains of intervals is given there.
What we add is an explicit, elementary order anti-isomorphism between
antichains of intervals and Dyck paths, passing through the staircase
functions, with complete proofs, together with
the translation of the rank of an antichain (computed in~\cite{BoVLAFI})
into the area of a Dyck path, which makes it possible to transfer level sets
and their asymptotics back and forth. We are not aware of previous work
connecting minimal-interval semantics with the combinatorics of Dyck paths.

\paragraph{Area and unimodality.} There is a bijection between Dyck paths of
semilength $m$ and Young diagrams fitting inside the staircase shape
$(m-1,m-2,\dots,1)$ that maps the area of a path to the difference between
the area of the staircase and the area of the diagram~\cite{FerUDP}; in
Section~\ref{sec:ccb} we obtained the same correspondence from the \CCB
lattice.

Both the unimodality of the Carlitz--Riordan polynomials and the Sperner
property of the Stanley lattice were first considered, to the best of our
knowledge, by Winston and Kleitman~\cite[Section~12]{WiKANTSS}; the
rank-unimodality of the Dyck lattices was later called an ``intriguing open
question'' by Bonin, Shapiro and Simion~\cite{BoSSQASNACSLP}, as reported
in~\cite{FePLLP}. Winston and Kleitman
observed, by computer, that the polynomials are unimodal up to semilength
$21$, and conjectured that they are unimodal in general (their Conjecture~2);
moreover, they considered the lattice of walks ordered by the relation of
lying above at each step, observed that the Dyck paths form an order ideal in
it, and wrote that ``it is not clear whether or not the smaller lattice has
the Sperner property''. Stanton~\cite{StaUYL}, who notes that the generating function of the partitions in a staircase was already considered by Carlitz~\cite{CarSPBN}, exhibited several infinite
families of partitions $\lambda$ for which the partitions fitting inside
$\lambda$ are not rank-unimodal, and conjectured (without citing~\cite{WiKANTSS}) that the
staircase shape $(n,n-1,\dots,1)$ is unimodal, verifying the conjecture for
$n\leq22$ (that is, for Dyck paths of semilength at most $23$); he also
reports that he and Fishel verified the Sperner property of the lattice of
partitions fitting inside any partition of size at most
$26$~\cite[p.~52]{StaUYL}, which includes the staircases $(n,n-1,\dots,1)$
with $n\leq6$. Ferrari~\cite{FerUDP}, who attributes the first appearance of the
conjecture to Stanton, describes its later restatements for the area of Dyck
paths and as rank-unimodality of the Stanley lattice, and proposes an approach
based on a partition of the Stanley lattice into saturated chains, without
settling the conjecture; we are not aware of a proof. Whether $\E_n$ is
Sperner, that is, whether its width equals the size $\lambda_n$ of its largest
level set, was asked again in~\cite{BoVLAFI} (where it was verified for
$n\leq 11$); our results show that it is the same question as the one of
Winston and Kleitman, whose answer is positive for $n\leq 20$ (Section~\ref{sec:partitions}). The literature on the largest coefficient of the Carlitz--Riordan
polynomial, which we use for the asymptotics of $\lambda_n$, was discussed in
Section~\ref{sec:asymptotics}.

\paragraph{Dimension and convex subsets.} The dimension of a poset is the
classical parameter defined by Dushnik and Miller~\cite{DuMPOS}, who also
proved that a poset has dimension at most two if and only if it can be
represented by a family of intervals of a linearly ordered set under
inclusion~\cite[Theorem~3.61]{DuMPOS}. In particular, the poset of nonempty
intervals of a chain has dimension at most two, in accordance with
Proposition~\ref{prop:dimsweep}, as the sweep dimension of a chain with at
least two elements is two. The intervals of a poset in
our sense are its \emph{convex} subsets, whose lattice $\mathrm{Co}(P)$
(ordered by inclusion) was studied by Birkhoff and Bennett~\cite{BiBCLP}. The
family of convex subsets of a finite poset is a \emph{convex geometry} (the
\emph{order convex alignment} of~\cite[Example~II]{EdJTCG}). Edelman and
Jamison~\cite{EdJTCG} introduced the \emph{convex dimension} $\mathrm{cdim}$
of a convex geometry, the minimum number of maximal chains of convex sets such
that every convex set is an intersection of members of the chains, and proved
that it is at least the size of the largest independent
set~\cite[Theorem~5.3]{EdJTCG}. Edelman and Saks~\cite[Theorem~4.1]{EdSCRCDCG}
proved that the convex dimension is equal to the width of the poset of
meet-irreducible convex sets, and observed that for every finite lattice $L$
the width of the poset of meet-irreducible elements is an upper bound, and the
height of the largest Boolean sublattice is a lower bound, for the
dimension of $L$~\cite[Proposition~4.4]{EdSCRCDCG}. Knauer
and Trotter~\cite{KnTCDCG} have recently compared several dimension parameters
of convex geometries, seen as posets under inclusion. In particular, they
record the inequality $\mathrm{cdim}\geq\dim$~\cite[Proposition~3.9]{KnTCDCG}
(which is also~\cite[Theorem~4.5]{EdSCRCDCG}), recall from Figure~1
of~\cite{EdSCRCDCG} that it can be strict~\cite[Section~4]{KnTCDCG}, and
show that the two parameters can be separated: there are convex geometries of
dimension three and arbitrarily large convex
dimension~\cite[Theorem~1.1]{KnTCDCG}.

Since the dimension of a poset and of its dual coincide,
Proposition~\ref{prop:dimsweep} states that
$\sweepdim(P,\sqleq)=\dim(\mathrm{Co}(P))$, so the sweep dimension is the
Dushnik--Miller dimension of a specific convex geometry, and the bounds of
Proposition~\ref{prop:sweepprops} can be recovered from the bounds of
Edelman and Saks. Indeed, the height of the largest Boolean sublattice of a
convex geometry is the maximum number of extreme points of a convex
set~\cite[Theorem~4.6]{EdSCRCDCG}, and the extreme points of a convex subset
$K$ of $P$ are its minimal and maximal elements~\cite[Example~II]{EdJTCG}. If
$S\subseteq P$ contains no $3$-element chain and $K$ is the smallest convex set
containing $S$, every $s\in S$ is minimal or maximal in $K$ (otherwise
$k_1\sqless s\sqless k_2$ for some $k_1,k_2\in K$, and since every element of
$K$ lies between two elements of $S$, we would find $s_1\sqless s\sqless s_2$
with $s_1,s_2\in S$), whence the lower bound
$\theta(P,\sqleq)\leq\sweepdim(P,\sqleq)$. For the upper bound, in a convex geometry the meet-irreducible
convex sets are the \emph{copoints}, that is, the maximal convex sets avoiding
some point $y$ (one direction is~\cite[Theorem~4.2]{EdJTCG}; the other is
immediate, as covers add a single point); a convex set avoiding $y$ cannot
contain both an element above $y$ and an element below $y$, so a copoint
attached at $y$ is either $P\setminus\upset y$ or $P\setminus\downset y$. As
$y$ varies, these two families are isomorphic to $(P,\sqleq)$ and to its dual,
respectively, so their union has width at most $2\cdot\width(P,\sqleq)$, which
yields
$\sweepdim(P,\sqleq)\leq\mathrm{cdim}(\mathrm{Co}(P))\leq2\cdot\width(P,\sqleq)$
by~\cite[Theorem~4.1 and Proposition~4.4]{EdSCRCDCG}.
We nonetheless give direct proofs, which in the case of the upper bound produce
an explicit sweep system. The inequality
$\sweepdim(P,\sqleq)\leq\mathrm{cdim}(\mathrm{Co}(P))$ can be strict: the
convex geometry of Figure~1 of~\cite{EdSCRCDCG}, which has convex dimension
four and dimension three (see also~\cite[Section~4]{KnTCDCG}), is the family of convex subsets of the four-element
poset $2\times2$, whose sweep dimension is thus three.

Edelman and Jamison also defined the \emph{interval dimension} of a poset
$(P,\sqleq)$ as the minimum number of linear orders $E_1$,~\dots,~$E_k$ on $P$
such that the convex subsets of $P$ are exactly the intersections
$I_1\cap\dots\cap I_k$, where each $I_i$ is an interval of
$E_i$~\cite[Section~5]{EdJTCG}; the definition makes sense because, for
every convex set $C$ and every $x\notin C$, there is a linear extension of
$P$ in which $x$ does not lie between two elements of~$C$~\cite[Lemma~5.8]{EdJTCG}.
This parameter is closely related to the sweep
dimension: since an interval of a linear order is the intersection of a final
segment of the order and of a final segment of its reverse, replacing each
$E_i$ with $E_i$ and its reverse yields a sweep system, so the sweep dimension
is at most twice the interval dimension. The two parameters are nonetheless
different: for an antichain with $n$ elements the interval dimension is
bounded by $\lfloor n/2\rfloor+1$~\cite{EdJTCG}, whereas the sweep dimension is $n$. We are not aware of previous work that specializes the Dushnik--Miller
dimension of a convex geometry to the poset of convex subsets of a poset, that
characterizes it in terms of sweeps or of orthants, or that studies the
staircase representation of antichains of intervals
(Theorem~\ref{thm:region}).

\section{Conclusions and open problems}
\label{sec:conclusions}

We have shown that the \CCB lattice is (up to its top element) anti-isomorphic to the Stanley lattice of Dyck words, and that under this correspondence the rank of an antichain is determined by the area of the corresponding Dyck word; this observation has allowed us to transfer to $\E_n$ the enumerative results available for Dyck words, and to obtain the asymptotics of its largest level set. For an arbitrary finite poset, we have introduced the notion of sweep dimension, and we have shown that it governs a representation of antichains of intervals as staircase regions; the sweep dimension coincides with the order dimension of the poset of intervals, and it is essentially determined by the width, up to a factor of two.

\paragraph{Consequences for retrieval.} Some consequences of the results of Section~\ref{sec:ccb} for minimal-interval semantics are worth stating. First, the Dyck word $\omega(A)$ is a balanced-parentheses sequence of length $2(n+1)$ that determines the antichain $A$ of minimal intervals over a document of $n$ positions: since $|\E_n|=C_{n+1}+1$, this encoding uses $2(n+1)$ bits, which is optimal up to lower-order terms for a worst-case answer, and balanced-parentheses sequences are a standard target of succinct data structures~\cite{JacSSTG,MuRSRBPST} (for sparse answers, of course, an explicit list of intervals can be shorter). Second, $\rho_A(\ell)$ is the smallest right end of a witness starting at or after $\ell$, which is the primitive that a lazy evaluation of the operators of minimal-interval semantics needs~\cite{BoVEOLAMIS}; the staircase function $\rho_A$ and the Dyck word $\omega(A)$ are two representations of the same primitive. Third, the rank of $A$ is $\sum_{x\in[n]}(n-\rho_A(x))$, a measure of the informativeness of the answer~\cite{BoVLAFI}, increasing along the order, and the level sets of $\E_n$ group the answers with equal informativeness; the asymptotics of $\lambda_n$ then says how large the largest group of equally informative answers is. This last statement is purely combinatorial: it counts all antichains uniformly, whereas the answers to actual queries on actual documents are far from uniformly distributed.

Several questions remain open. 
\begin{itemize}
\item \emph{The value of $\lambda_n$.} We have obtained the asymptotics of $\lambda_n$, and Table~\ref{tab:lambda} shows that the convergence is slow; a closed formula, or an asymptotic expansion with a second term, would be desirable. Numerically, the relative error of the first-order estimate decays like $1/n$ (see Section~\ref{sec:asymptotics}).
\item \emph{The sweep dimension of specific posets.} We have computed the sweep dimension of chains, antichains and complete bipartite posets. 
It would be interesting to determine it for other natural classes, for instance trees (Figure~\ref{fig:orthant} shows that the sweep dimension of a tree 
is not governed by the intervals of the form $[a\..b]$), products of chains, or Boolean lattices, and to characterize the posets 
for which the sweep dimension attains the lower bound $\max(2,\width)$ or the upper bound $2\cdot\width$.
\item \emph{Which staircase regions arise?} For a chain, the map $A\mapsto R(A)$ is a bijection onto the staircase regions contained in a triangle. For a general poset and a given sweep system, a characterization of the regions of the form $R(A)$ would provide an analogue of the Dyck-word representation for arbitrary posets, and possibly a way to enumerate $\E(P,\sqleq)$. Since the regions depend on the chosen sweep system, it would also be interesting to find a canonical choice, or to understand how the family of regions varies among the realizers of $(\I(P,\sqleq),\supseteq)$ of minimum size.
\end{itemize}
\backmatter

\bmhead{Funding}
The authors did not receive support from any organization for the submitted work.

\bmhead{Conflict of interest}
The authors declare no competing interests.

\bmhead{Data availability}
Data sharing is not applicable to this article as no datasets were generated or analysed during the current study.

\bibliography{Previous/biblio,extra}

\end{document}